\documentclass[11pt]{amsart}

\usepackage{amsmath, amsthm, amssymb, amscd, amsfonts}
\usepackage{fontsmpl}
\usepackage{fontsmpl,layout}
\usepackage[all]{xy}

\newtheorem{mtheorem}{Theorem}
\newtheorem{lema}{Lemma}[section]
\newtheorem{theorem}[lema]{Theorem}
\newtheorem{cor}[lema]{Corollary}

\newtheorem{prop}[lema]{Proposition}

\theoremstyle{definition}
\newtheorem{definition}[lema]{Definition}

\theoremstyle{remark}
\newtheorem{obs}[lema]{Remark}
\newtheorem{rmk}[lema]{Remarks}

\newcommand{\quot}{\mathcal{QUOT}}

\newcommand\id{\operatorname{id}}

\newcommand\gr{\operatorname{gr}}
\newcommand\co{\operatorname{co}}

\newcommand\rk{\operatorname{rk}}

\newcommand\Hom{\operatorname{Hom}}
\newcommand\Alg{\operatorname{Alg}}

\newcommand\Ker{\operatorname{Ker}}
\newcommand\Img{\operatorname{Im}}

\newcommand{\eps}{\varepsilon}
\newcommand{\ot}{\otimes}
\newcommand{\com}{\Delta}

\newcommand\Res{\operatorname{Res}}
\newcommand\res{\operatorname{res}}
\newcommand{\g}{\mathfrak g}

\newcommand{\sop}{\operatorname{Supp}}

\newcommand\toba{{\mathfrak B }}

\newcommand{\ku}{\mathbb C}

\newcommand{\D}{{\mathcal D}}

\newcommand{\Pc}{{\mathcal P}}
\newcommand{\Oc}{{\mathcal O}}

\newcommand{\yduo}{{}^{U_{0}}_{U_{0}}\mathcal{YD}}
\newcommand\Ent{\operatorname{Ent}}
\newcommand\Fr{\operatorname{Fr}}
\newcommand\Lie{\operatorname{Lie}}
\renewcommand{\_}[1]{_{\left( #1 \right)}}

\def\NN{\mathbb{N}}
\def\ZZ {\mathbb{Z}}
\def\QQ{\mathbb{Q}}
\def\CC{\mathbb{C}}

\def\AA{\mathbb{A}}
\def\JJ{\mathcal{J}}
\def\II{\mathcal{I}}
\def\OO{\mathcal{O}}
\def\SS{\mathcal{S}}

\def\z{\mathbf{Z}}
\def\x{\mathbf{X}}
\def\y{\mathbf{Y}}

\def\D{\mathcal{D}}

\def\J{\mathfrak{J}}

\def\lieg{\mathfrak{g}}
\def\lieh{\mathfrak{h}}
\def\liek{\mathfrak{k}}

\def\liek{\mathfrak{k}}

\def\QEp{{\check{U}}_{q}(\lieg)}
\def\QEpe{{\check{U}}_{\epsilon}(\lieg)}

\def\Rq{R_{q}[G]}
\def\Ga{\Gamma}

\def\Gli{\Ga(\lieg)}
\def\Glie{\Ga_{\epsilon}(\lieg)}
\def\Oeq{\OO_{\epsilon}(G)_{\QQ(\epsilon)}}
\def\Oe{\OO_{\epsilon}(G)}

\def\Oesln{\OO_{\epsilon}(SL_{N})}
\def\Oec{\OO_{\epsilon}(G)_{\CC}}
\def\Oek{\OO_{\epsilon}(G)_{k}}

\def\qe{{\bf u}_{\epsilon}(\lieg)}

\def\liel{\mathfrak{l}}
\def\lgot{\mathfrak{l}}
\def\agot{\mathfrak{a}}
\def\bgot{\mathfrak{b}}
\def\QEbmas{{\check{U}}_{\epsilon}(\bgot_{+})}
\def\QEbmen{{\check{U}}_{\epsilon}(\bgot_{-})}
\def\Ege{U(\lieg)_{\QQ(\epsilon)}}
\def\Ele{U(\liel)_{\QQ(\epsilon)}}
\def\Gl{\Ga(\liel)}

\def\Glil{\Ga_{\epsilon}(\liel)}
\def\Oeql{\OO_{\epsilon}(L)_{\QQ(\epsilon)}}
\def\Oel{\OO_{\epsilon}(L)}
\def\OG{\OO(G)_{\QQ(\epsilon)}}
\def\OL{\OO(L)_{\QQ(\epsilon)}}

\def\Ol{\Oc_{\epsilon}(L)}
\def\Olp{\Oc_{\epsilon}(L')}
\def\ql{{\bf u}_{\epsilon}(\lgot)}
\def\qlp{{\bf u}_{\epsilon}(\lgot')}

\def\pf{\begin{proof}}
\def\epf{\end{proof}}
\theoremstyle{plain}

\begin{document}




\title[quantum subgroups]
{Quantum subgroups of a simple quantum group at roots of 1}

\author[n. andruskiewitsch and g. a. garc\'\i a]
{Nicol\' as Andruskiewitsch  \and Gast\' on Andr\' es Garc\'\i a}

\address{FaMAF-CIEM, Universidad Nacional de C\'ordoba
Medina Allende s/n, Ciudad Universitaria, 5000 C\' ordoba, Rep\'
ublica Argentina.} \email{andrus@mate.uncor.edu,
ggarcia@mate.uncor.edu}

\thanks{Results
in this paper are part of the Ph.D. thesis of G. A. G., written
under the advise of N. A. The work was partially supported by
CONICET, ANPCyT, Secyt (UNC) and TWAS}

\date{June 22, 2007.}

\begin{abstract}
Let $G$ be a connected, simply connected, simple complex algebraic
group and let $\epsilon$ be a primitive $\ell$-th root of 1, $\ell$
odd and $3 \nmid \ell$ if $G$ is of type $G_{2}$. We determine all
Hopf algebra quotients of the quantized coordinate algebra $\Oe$.
\end{abstract}

\maketitle



\section{Introduction and preliminaries}

\subsection{Introduction}
The purpose of this paper is to determine all quantum subgroups of
a quantum group at a root of one, or in equivalent terms, to
determine all Hopf algebra quotients of a quantized coordinate
algebra at a root of one (over the complex numbers). This problem
was first considered by P. Podle\'s \cite{podles} for quantum
$SU(2)$ and $SO(3)$. The characterization of all
\emph{finite-dimensional} Hopf algebra quotients of the quantized
coordinate algebra $\Oc_q(SL_N)$ was obtained by Eric M\"uller
\cite{Mu2}. M\"uller's approach is via explicit computations with
matrix coefficients; this strategy does not apply to more general
simple groups.

\smallbreak The present work can be viewed as a continuation of the
long tradition of studying subgroups of a simple algebraic group. In
fact, our main theorem assumes the knowledge of such subgroups, see
Definition \ref{def:subgroupdatum}. Besides its intrinsical
mathematical interest, our result would have implications in quantum
harmonic analysis-- see for example \cite{letzter}-- and in the
study of module categories over the tensor category of comodules
over the Hopf algebra $\Oe$-- in the sense of \cite{eo}.

\smallbreak An outcome of our main theorem is the construction of
many new examples of finite-dimensional Hopf algebras. At the
present time, all examples of finite-dimensional Hopf algebras, we
are aware of, are:
\begin{itemize}
    \item group algebras of finite groups,
   \item small quantum groups introduced by Lusztig \cite{L1, L},
and variations thereof \cite{AS-05},
\item other pointed Hopf algebras with abelian group
arising from the Nichols algebras discovered in \cite{G1, He},
\item a few examples of pointed Hopf algebras with non-abelian group
\cite{MS, G2},
    \item combinations of the preceding via some standard
    operations (duals, twisting, Hopf subalgebras and quotients, extensions).
\end{itemize}

How to build examples of Hopf algebras via extensions of a group
algebra by a dual group algebra is well understood-- see for
instance \cite{ma-newdir}. Out of this, extensions can in
principle be constructed by means of weak actions and coactions,
and pairs of compatible 2-cocycles. However, very few explicit
examples were presented in this way, to our knowledge no one in
finite dimension, except for the trivial tensor product of two
Hopf algebras. Our examples are indeed nontrivial extensions of
finite quantum groups by finite groups, but it is not clear how
they could be explicitly presented through actions, coactions and
cocycles. A natural subsequent question is when the new examples
of Hopf algebras are isomorphic with each other; this will be
addressed in (the forthcoming new version of) \cite{AG}.

Furthermore, a result of \c Stefan \cite[Thm. 1.5]{stefan} says that
a non-semi\-simple finite-dimensional Hopf algebra generated by a
simple 4-dimensional coalgebra stable by the antipode is a quotient
of the quantized coordinate algebra of $SL(2)$ at a root of one. It
is tempting to suggest that finite-dimensional quotients of more
general quantized coordinate algebras might play a prominent role in
the classification of Hopf algebras.

\medbreak We notice that a different problem is sometimes referred
to with a similar name: this is the classification of
indecomposable module categories over fusion categories arising in
conformal field theory, \emph{e. g.} from the representation
theory of finite quantum groups at roots of one. See
\cite{ocneanu, ko}. There is no evident relation between these two
problems.

\subsection{Statement of the main result} Let $\g$ be the
Lie algebra of $G$, $\lieh\subseteq \g$ a fixed Cartan subalgebra,
$\Pi=\{\alpha_{1},\ldots,\alpha_{n} \}$ a basis of the root system
$\Phi = \Phi(\g,\lieh)$ of $\g$ with respect to $\lieh$ and $n =
\rk \g$.
\begin{definition}\label{def:subgroupdatum}
A \emph{subgroup datum} is a collection $\D= (I_{+}, I_{-}, N, \Ga,
\sigma, \delta)$ where
\begin{itemize}
\item[$\bullet$] $I_{+}\subseteq \Pi$ and $I_{-} \subseteq -\Pi$.
Let $\Psi_{\pm} = \{\alpha \in \Phi:\ \sop \alpha \subseteq
I_{\pm}\}$, $\lgot_{\pm} = \sum_{\alpha \in \Psi_{\pm}} \g_{\alpha}$
and $\lgot = \lgot_{+} \oplus \lieh \oplus \lgot_{-}$; $\lgot$ is an
algebraic Lie subalgebra  of $\g$. Let $L$ be the connected Lie
subgroup of $G$ with $\Lie (L) = \lgot$.

\medbreak \item[$\bullet$] $N$ is a subgroup of
$\widehat{\mathbb{T}_{I^{c}}}$, see Remark \ref{rmk:subgruposigma}
below.

\medbreak \item[$\bullet$] $\Ga$ is an algebraic group.

\medbreak \item[$\bullet$] $\sigma: \Ga \to L$ is an injective
homomorphism of algebraic groups.

\medbreak \item[$\bullet$] $\delta: N \to \widehat{\Ga}$ is a
group homomorphism.
\end{itemize}
\end{definition}

\noindent If $\Ga$ is finite, we call $\D$ a \emph{finite subgroup
datum}. We parameterize with injective group homomorphisms rather
than group inclusions for a better description of the isomorphism
classes \cite{AG}. An equivalence relation among subgroup data is
defined in Subsection \ref{sec:equivalence}.

\bigbreak Our main result is the following.

\begin{mtheorem}\label{teo:biyeccion}
There is a bijection between
\begin{enumerate}
\item[$(a)$] Hopf algebra quotients $q: \Oe \to A$.
\item[$(b)$] Subgroup data up to equivalence.
\end{enumerate}
\end{mtheorem}

In Section \ref{section:q-subgroups}, we carry out the construction
of a quotient $A_{\D}$ of $\Oe$ starting from a subgroup datum $\D$,
see Theorem \ref{teo:constrfinal}. In Subsection
\ref{sec:equivalence}, we study the lattice of quotients $A_{\D}$.
In Section \ref{sec:demmainthm}, we attach a subgroup datum $\D$ to
an arbitrary Hopf algebra quotient $A$ and prove that $A_{\D}\simeq
A$ as quotients of $\Oe$. This concludes the proof of the Theorem
\ref{teo:biyeccion}. As an immediate corollary of Theorem
\ref{teo:biyeccion}, we get
\begin{mtheorem}\label{teo:biyeccionfin}
There is a bijection between
\begin{enumerate}
\item[$(a)$] Hopf algebra quotients $q: \Oe \to A$ such that
$\dim A < \infty$.
\item[$(b)$] Finite subgroup data up to equivalence.
\end{enumerate}
\end{mtheorem}

Theorem \ref{teo:biyeccionfin} generalizes the main result of
\cite{Mu2}.

\subsection{Conventions}\label{sec:prelim}
Let $C = (a_{ij})_{1\leq i,j\leq n}$ be the Cartan matrix of $\g$
and suppose that $\g$ is generated by the elements $\{h_{i},\
e_{i},\ f_{i}\vert\ 1\leq i\leq n\}$ subject to the
Chevalley-Serre relations. Let $Q = \ZZ \Phi = \bigoplus_{i=1}^{n}
\ZZ \alpha_{i}$ be the root lattice, $\varpi_{1},\ldots
,\varpi_{n}$ the fundamental weights, $P = \bigoplus_{i=1}^{n}
\ZZ\varpi_{i}$ the weight lattice and $W$ the Weyl group. Let
$P_+$ be the cone of dominant weights and $Q_+ = P_+ \cap P$. Let
$(-,-)$ be the positive definite symmetric bilinear form on
$\lieh^{*}$ induced by the Killing form of $\lieg$. Let $d_{i} =
\frac{(\alpha_{i},\alpha_{i})}{2} \in \{1,2,3\}$.

\smallbreak For $t,\ m \in \NN_{0}$, $q \in \CC$ and $u \in
\QQ(q)\smallsetminus \{0,\pm 1\}$ we denote:
\begin{align*} [t]_{u}: & = \frac{u^{t}-u^{-t}}{u-u^{-1}},&
[t]_{u}! & := [t]_{u}[t-1]_{u}\cdots [1]_{u}, & \left[\begin{matrix}m\\
t\end{matrix}\right]_{u} & := \frac{[m]_{u}!}{[t]_{u}!
[m-t]_{u}!},\\
(t)_{u} & : = \frac{u^{t} -1}{u -1}, & (t)_{u}!  & :=
(t)_{u}(t-1)_{u}\cdots (1)_{u},&
\left(\begin{matrix}m\\
t\end{matrix}\right)_{u}  & := \frac{(m)_{u}!}{(t)_{u}! (m-t)_{u}!}.
\end{align*}

\subsection{Definitions}
In this subsection we recall the definition of the quantized
coordinate algebra of $G$. Let $R = \QQ[q,q^{-1}]$, $q$ an
indeterminate. If $p_{\ell}(q) \in R$ denotes the $\ell$-th
cyclotomic polynomial, then $R/[p_{\ell}(q)R] \simeq \QQ(\epsilon)$.

\begin{definition} The {\it simply connected}
quantized enveloping algebra $\QEp$ of $\lieg$ is the
$\QQ(q)$-algebra with generators $\{K_{\lambda}\vert\ \lambda \in
P\}$, $E_{1},\ldots, E_{n}$ and $F_{1},\ldots, F_{n}$, satisfying
the following relations for $\lambda,\ \mu \in P$ and $1\leq i,j\leq
n$:

\begin{align*}
K_{0} = 1, &\qquad K_{\lambda}K_{\mu} = K_{\lambda+\mu}, \\
K_{\lambda}E_{j}K_{-\lambda}  = q^{(\lambda,\alpha_{j})}E_{j}, &
\qquad K_{\lambda}F_{j}K_{-\lambda}
=   q^{-(\lambda, \alpha_{j})}F_{j},\\
E_{i}F_{j} - F_{j}E_{i} & =
\delta_{ij}\frac{K_{\alpha_{i}}-K_{\alpha_{i}}^{-1}}{q_{i} -
q_{i}^{-1}},\\
\sum_{l=0}^{1-a_{ij}} (-1)^{l}&
{\left[ \begin{smallmatrix} 1-a_{ij}\\
l\end{smallmatrix} \right]_{q_{i}}} E_{i}^{1-a_{ij}-l}E_{j}E_{i}^{l}
= 0 \qquad (i\neq j),\\
\sum_{l=0}^{1-a_{ij}} (-1)^{l}& {\left[ \begin{smallmatrix} 1-a_{ij}\\
l \end{smallmatrix}\right]_{q_{i}}} F_{i}^{1-a_{ij}-l}F_{j}F_{i}^{l}
= 0 \qquad  (i\neq j).
\end{align*}
\end{definition}

\begin{definition} \cite[Section 3.4]{DL}
Let $q_{i} = q^{d_{i}}$, $1\leq i\leq n$. The algebra $\Gli$ is the
$R$-subalgebra of $\QEp$ generated by the elements

\begin{align*}
K_{\alpha_{i}}^{-1} & & (1\leq i\leq n),\\
\left(\begin{array}{c}K_{\alpha_{i}}; 0\\
t\end{array}\right)&:=
\prod_{s=1}^{t}\left(\frac{K_{\alpha_{i}}
q_{i}^{-s+1}-1}{q_{i}^{s}-1}\right)
& (t\geq 1,\ 1\leq i\leq n),
\\
E_{i}^{(t)} & := \frac{E_{i}^{t}}{[t]_{q_{i}}!} &
 (t\geq 1,\ 1\leq i\leq n),\\
F_{i}^{(t)} &:= \frac{F_{i}^{t}}{[t]_{q_{i}}!} &
 (t\geq 1,\ 1\leq i\leq n).
\end{align*}
\end{definition}

Let $\mathcal{C}$ be the strictly full subcategory of $\Gli$-mod
whose objects are $\Gli$-modules $M$ such that $M$ is a free
$R$-module of finite rank and the operators
$K_{\alpha_{i}}$ and $\left(\begin{smallmatrix}K_{\alpha_{i}}; 0\\
t\end{smallmatrix}\right) $ are diagonalizable with
eigenvalues $q_{i}^{m}$ and $\left(\begin{smallmatrix}m\\
t\end{smallmatrix}\right)_{q_{i}}$ respectively, for some $m \in
\NN$ and for all $1\leq i\leq n$.

\begin{definition}\cite[Section 4.1]{DL}
Let $\Rq$ denote the $R$-submodule of $\Hom_{R}(\Gli,R)$ spanned by
the coordinate functions $t^{j}_{i}$ of representations $M$ from
$\mathcal{C}$: $<g,t^{j}_{i}> = <g\cdot m_{i}, m^{j} >$, where
$(m_{i})$ is an $R$-basis of $M$, $(m^{j})$ is the dual basis of the
dual module and $g \in \Gli$. Since the subcategory $\mathcal{C}$ is
a tensor one, $\Rq$ is a Hopf algebra.
\end{definition}

\begin{definition}\cite[Section 6]{DL}
The algebra $\Rq / [p_{\ell}(q)\Rq]$ is denoted by $\Oeq$ and is
called the {\it quantized coordinate algebra} of $G$ over
$\QQ(\epsilon)$ at the root of unity $\epsilon$. In the same way
as for $\Oeq$, we can form the $\QQ(\epsilon)$-Hopf algebra $\Glie
:= \Gli / [p_{\ell}(q)\Gli]$.
\end{definition}

We now relate the Hopf algebras $\Oeq$ and $\Glie$.

\begin{definition}
A {\it Hopf pairing} between two Hopf algebras $U$ and $H$ over a
ring $\mathcal{R}$ is a bilinear form $(-,-): H \times U \to
\mathcal{R}$ such that, for all $u,\ v \in U$ and $f,\ h \in H$,
\begin{align*}
& (i)\qquad  (h,uv) = (h_{(1)},u)(h_{(2)},v);\qquad & (iii)&\qquad
(1,u) =
\eps(u);\\
&(ii)\qquad (fh,u) = (f,u_{(1)})(h,u_{(2)}); \qquad & (iv)& \qquad
(h,1) = \eps(h).
\end{align*}
\end{definition}

\noindent It follows that $(h,\SS(u)) = (\SS(h),u)$, for all $u \in
U$, $h \in H$. Given a Hopf pairing, one has Hopf algebra maps $U
\to H^{\circ}$ and $H \to U^{\circ}$, where $H^{\circ}$ and
$U^{\circ}$ are the Sweedler duals. The pairing is called {\it
perfect} if these maps are injections.

\begin{prop}\label{perfectpairing}
\cite[4.1 and 6.1]{DL} There exists a perfect Hopf pairing $\Rq
\otimes_{R} \Gli \to R$, which induces a perfect Hopf pairing $\Oeq
\otimes_{\QQ(\epsilon)} \Glie \to \QQ(\epsilon)$. In particular,
$\Oeq \subseteq \Glie^{\circ}$ and $\Glie \subseteq \Oeq^{\circ}$.
\qed
\end{prop}

If $k$ is any field containing $\QQ(\epsilon)$, we denote $\Oek:=
\Oeq\otimes_{\QQ(\epsilon)} k$. When $k=\CC$ we simply write $\Oe$
for $\Oec$. The following two results imply by \cite[Prop.
3.4.3]{Mo} that $\Oe$ is a central extension of $\OO(G)$ by a
finite-dimensional Hopf algebra.

\begin{theorem}\label{OcentralenOe}
\begin{enumerate}
\item[$(a)$] \cite[Prop. 6.4]{DL}
$\Oe$ contains a central Hopf subalgebra isomorphic to the
coordinate algebra $\OO(G)$ of $G$.
\item[$(b)$] \cite[III.7.11]{BG} $\Oe$ is a free
$\OO(G)$-module of rank $\ell^{\dim G}$. \qed
\end{enumerate}
\end{theorem}

We end this section by spelling out explicitly the quotient of $\Oe$
by its central Hopf subalgebra $\OO(G)$.

\smallbreak Let $\overline{\Oe}= \Oe / [\OO(G)^{+}\Oe]$ and denote
by $\pi: \Oe \to \overline{\Oe}$ the quotient map. By Theorem
\ref{OcentralenOe} and \cite[Prop. 3.4.3]{Mo}, $\overline{\Oe}$ is a
Hopf algebra of dimension $\ell^{\dim G}$ which fits into the exact
sequence
$$ 1\to \OO(G) \to \Oe \to \overline{\Oe} \to 1.$$

\noindent Let $\qe$ be the {\it Frobenius-Lusztig kernel} of
$\lieg$ at $\epsilon$; that is, the Hopf subalgebra of $\Glie$
generated by the elements $E_{i}$, $F_{i}$ and $K_{\alpha_{i}}$
for $1\leq i\leq n$. See \cite{BG} for details. We denote by
\begin{equation}\label{eqn:torofinito} \mathbb{T}
:=\{K_{\alpha_{1}},\ldots ,K_{\alpha_{n}}\} = G(\qe)
\end{equation}

\noindent the ``finite torus'' of group-like elements of $\qe$.

\begin{theorem}\label{isoOeqe} \cite[III.7.10]{BG}
$\overline{\Oe} \simeq \qe^{*}$ as Hopf algebras. \qed
\end{theorem}

Summarizing, the quantized coordinate algebra $\Oe$ of $G$ at
$\epsilon$ fits into the central exact sequence
\begin{equation}\label{extOOeqe}
 1\to \OO(G) \xrightarrow{\iota} \Oe \xrightarrow{\pi} \qe^{*} \to 1.
\end{equation}

\noindent We shall need  the following technical lemma.

\begin{lema}\label{lema:quotient}
There exists a surjective algebra map $\varphi: \Glie \to \qe$
such that $\varphi|_{\qe} = \id$.
\end{lema}

\pf Since $\Glie = \Gli / [p_{\ell}(q)\Gli]$, we may define
$\varphi$ as a map from $\Gli$ such that $\varphi(q) = \epsilon$.
Let $\varphi$ be the unique algebra map which takes the following
values on the generators:

\begin{align*}
\varphi(E_{i}^{(m)}) & =  \begin{cases}{\begin{matrix}E_{i}^{(m)} &
\mbox{if }
1\leq m < \ell\\
0& \mbox{otherwise,}\end{matrix}}\end{cases} \\
\varphi(F_{i}^{(m)}) & =
\begin{cases}\begin{matrix}F_{i}^{(m)} & \mbox{if } 1\leq m < \ell
\\
0& \mbox{otherwise,}\end{matrix}\end{cases}\\
\varphi(\left(\begin{smallmatrix} K_{\alpha_{i}};0\\
m\end{smallmatrix}\right)) & =
\begin{cases}\begin{matrix} \left(\begin{smallmatrix} K_{\alpha_{i}};0\\
m\end{smallmatrix}\right)& \mbox{if } 1\leq m < \ell \\
0& \mbox{otherwise,}\end{matrix}\end{cases}\\
\varphi(K_{\alpha_{i}}^{-1}) & = K_{\alpha_{i}}^{\ell
-1},\qquad\qquad \varphi(q) = \epsilon,
\end{align*}

\noindent for all $1\leq i\leq n$. Since $\varphi$ is the identity
on the generators of $\qe$ and $E_{i}^{\ell} = 0 = F_{i}^{\ell}$,
$K_{\alpha_{i}}^{\ell} = 1$ on $\qe$, it follows from a direct
computation that $\varphi$ satisfies the relations given in
\cite[Section 3.4]{DL}, see \cite[4.1.17]{G} for details. Hence
$\varphi$ is a well-defined algebra map whose image is $\qe$. \epf

\subsection{Hopf subalgebras of a pointed Hopf algebra}
\label{section:hopf-subalg} We describe in this subsection Hopf
subalgebras of pointed Hopf algebras. Let $U$ be a Hopf algebra such
that the coradical $U_0$ is a Hopf subalgebra. Let $(U_n)_{n\ge 0}$
be the coradical filtration of $U$, set $U_{-1} = 0$, $\gr U(n) =
U_n/U_{n-1}$ and let $\gr U = \oplus_{n\ge 0}\gr U(n)$ be the
associated graded Hopf algebra. Let $\iota: U_0\to\gr U$ be the
canonical inclusion and let $\pi:\gr U\to U_0$ be the homogeneous
projection. Let $R= (\gr U)^{\co \pi}$ be the diagram of $U$; $R$ is
a graded braided Hopf algebra, that is, a Hopf algebra in the
category $\yduo$ of Yetter-Drinfeld modules over $U_{0}$. Its
coalgebra structure is given by $\com_{R}(r) =
\vartheta_R(r\_{1})\ot r\_{2}$, for all $r\in R$, where
$\vartheta_R: \gr U \to R$ is the map defined by
\begin{equation}\label{vartheta}\vartheta_R(a) = a\_{1}\iota\pi(\SS
a\_{2}),\qquad \forall\ a\in \gr U.\end{equation}

\noindent It can be easily shown that $\vartheta_R(rh) =
r\epsilon(h)$, $\vartheta_R(hr) = h\cdot r$ for $r\in R$, $h\in
U_0$. One has that $\gr U \simeq R\# U_0$, $R = \oplus_{n\ge
0}R(n)$, $R(0) \simeq \ku$ and $R(1) = \Pc(R)$. We say that $R$ is a
\emph{Nichols algebra} if $R$ is generated as algebra by $R(1)$. See
\cite{AS1} for more details.

To state the following result, we need to introduce some
terminology. Let $A$ be a Hopf algebra, $M$ a Yetter-Drinfeld module
over $A$ and $B$ a Hopf subalgebra of $A$. We say that a vector
subspace $N$ of $M$ is \emph{$B$-compatible} if
\begin{enumerate}
    \item[$(a)$] it is stable under the action of $B$, and
    \item[$(b)$] it bears a $B$-comodule structure inducing the
    coaction of $A$.
\end{enumerate}

In inaccurate but descriptive words, ``$N$ is a Yetter-Drinfeld
submodule over $B$" (although $M$ is not necessarily a
Yetter-Drinfeld module over $B$).

\begin{lema}\label{lema:hopf-subalg} Let $Y$ be a Hopf subalgebra of
$U$. Then the coradical $Y_0$ is a Hopf subalgebra and the diagram
$S$ of $Y$ is a braided Hopf subalgebra of $R$.

If $R = \toba(V)$ is a Nichols algebra with $\dim V < \infty$,
then $S$ is also a Nichols algebra. In this case, Hopf subalgebras
of $U$ are parameterized by pairs $(Y_0, W)$ where $Y_0$ is a Hopf
subalgebra of $U_0$ and $W\subset V = R(1)$ is $Y_0$-compatible.
\end{lema}

\pf The first claim follows since $Y_0 = Y\cap U_0$ and the
intersection of two Hopf subalgebras is a Hopf subalgebra. By
\cite[Lemma 5.2.12]{Mo}, the coradical filtration of $Y$ is given by
$Y_n = Y\cap U_n$; thus we have an injective homogeneous map of Hopf
algebras $\gamma: \gr Y \hookrightarrow \gr U$ inducing the
commutative diagram
$$
\xymatrix{\gr Y\ar@{^{(}->}^{\gamma}[0,2]\ar@{->}[1,0]^{\pi_Y}&
& \gr U\ar@{->}[1,0]^{\pi} \\
Y_0 \ar@{^{(}->}[0,2]& & U_0.}
$$
Thus $S = \{a\in \gr Y: (\id\ot\pi_Y)\Delta(a) = a\ot 1\}$ is a
subalgebra, and also a braided vector subspace, of $R$. Note that
$\gamma\vartheta_S = \vartheta_R\gamma$, \emph{cf.}
\eqref{vartheta}; thus $S$ is a subcoalgebra of $R$. Assume now
that $R\simeq \toba(V)$ is a Nichols algebra with $\dim V<
\infty$. Taking graded duals, we have a surjective map of graded
braided Hopf algebras $\wp:\toba(V^*)\to S^{\gr \text{dual}}$.
Since $\toba(V^*)$ and $S^{\gr \text{dual}}$ are pointed
irreducible coalgebras, by \cite[Thm. 9.1.4]{Sw}, $\wp$ maps the
coradical filtration of the first onto the coradical filtration of
the second; hence $\Pc(S^{\gr \text{dual}}) = S^{\gr
\text{dual}}(1)$ and \emph{a fortiori} $S$ is generated in degree
1, \emph{i.~e.} is a Nichols algebra. Furthermore, $Y$ is
determined by $Y_0$ and $S(1)$, the last being $Y_0$-compatible.
Conversely, if $Y_0$ is a Hopf subalgebra of $U_0$ and $W\subset
R(1)$ is $Y_0$-compatible, then choose $(y_i)_{i\in I}$ in $U_1$
such that the classes $(\overline{y}_i)_{i\in I}$ in $U_1/U_0$
generate $W\# 1$. Then the subalgebra $Y$ of $U$ generated by
$Y_0$ and $(y_i)_{i\in I}$ is a actually a Hopf subalgebra giving
rise to the pair $(Y_0, W)$. \epf

The lemma above also holds if $V$ is a locally finite braided vector
space.

\smallbreak Let us now turn to Hopf subalgebras of pointed Hopf
algebras. The notion of ``compatibility" for groups reads as
follows. Let $G$ be a group and $M$ a Yetter-Drinfeld module over
the group algebra $\ku [G]$. If $F$ is a subgroup of $G$, a vector
subspace $N$ of $M$ is $F$-\emph{compatible} if
\begin{enumerate}
    \item[$(a)$] it is stable under the action of $F$, and
    \item[$(b)$] it is a $\ku [G]$-subcomodule and
    $\sop N := \{g\in G: N^g\neq 0\}$ is contained in $F$.
\end{enumerate}

\begin{cor}\label{prop:hopf-subalg}
Let $U$ be a pointed Hopf algebra whose diagram $R$ is a
Nichols algebra. Then Hopf subalgebras of $U$ are parameterized by
pairs $(F, W)$ where $F$ is a subgroup of $G(U)$ and $W\subset R(1)$
is $F$-compatible. \qed
\end{cor}

The Corollary reads even nicer if $G(U)$ is abelian and $\dim R(1)^g
= 1$ for all $g\in \sop R(1)$. Indeed, Hopf subalgebras of $U$ are
parameterized in this case by pairs $(F, J)$ where $F$ is a subgroup
of $G(U)$ and $J\subset \sop R(1)$ is contained in $F$. We recover
in this way results from \cite{CM, MuI}.

\begin{cor}\label{cor:hopf-subalg-uepsilong}\cite[Thm. 6.3]{MuI}
The Hopf subalgebras of $\qe$ are parameterized by triples
$(\Sigma, I_+, I_-)$, where $\Sigma$ is a subgroup of $\mathbb{T}$
and $I_{+} \subseteq \Pi$, $I_{-} \subseteq -\Pi$ such that
$K_{\alpha_{i}} \in \Sigma$ if $\alpha_{i} \in I_+\cup
-I_{-}$.\qed
\end{cor}

\subsection{A five-lemma for extensions of Hopf algebras}
\label{section:five-lemma}

 The following general lemma was kindly
communicated to us by Akira Masuoka.

\begin{lema}\label{lema:iso-galois-extensions}
Let  $H$  be a bialgebra over an arbitrary commutative ring, and
let  $A$,  $A'$  be right $H$-Galois extensions over a common
algebra  $B$  of $H$-coinvariants. Assume that $A'$  is right
$B$-faithfully flat. Then any $H$-comodule algebra map $\theta : A
\to A'$ that is identical on $B$  is an isomorphism. \end{lema}

\pf Let $\beta : A \ot_{B} A \to A \ot H$, $\beta(x\ot y) = x
y_{(0)}\ot y_{(1)}$ and $\beta' : A' \ot_{B} A' \to A' \ot H$,
$\beta'(x'\ot y') = x' y'_{(0)}\ot y'_{(1)}$ be the corresponding
Galois maps, for $x,y\in A$, $x',y'\in A'$. Using the $A$-module
structure of $A'$ given by $\theta$, we can extend $\beta$ to an
isomorphism
$$\alpha: A'\ot_{B} A \simeq A'\ot_{A} A \ot_{B} A
\xrightarrow{\id\ot\beta}  A'\ot_{A} A \ot H \simeq A' \ot H.$$
Explicitly, $\alpha(a'\ot a) = a'\theta(a_{(0)})\ot a_{(1)}$ for
all $a' \in A'$, $a\in A$. Then $\alpha$ fits into the following
commutative diagram
$$\xymatrix{A'\ot_{B} A \ar[rr]^{\id\ot \theta}
\ar[dr]_{\alpha}^{\simeq}
& &  A'\ot_{B} A' \ar[dl]^{\beta}_{\simeq}\\
& A'\ot H & }$$ \noindent Hence $\id\ot\theta$ is an isomorphism;
since $A'$ is right $B$-faithfully flat, $\theta$ is an
isomorphism. \epf

The lemma applies to a commutative diagram of Hopf algebras
\begin{equation}\label{diag:five}
\xymatrix{1 \ar[r]^{} &   B \ar[r]^{{\iota}}\ar@{=}[d] & A
\ar[r]^{{\pi}}\ar@{->>}[d]^{\theta} & H\ar[r]^{}\ar@{=}[d] & 1\\
1 \ar[r]^{} &  B  \ar[r]^{\iota'} & A' \ar[r]^{\pi'}& H\ar[r]^{}&
1,}
\end{equation}
where the rows are exact sequences of Hopf algebras, in the sense of
\cite{AD}: $A^{\co \pi} = B$ and $\ker \pi = B^+A$; \emph{ditto} for
$A'$. If the top row is a cleft exact sequence, then $\theta$ is an
isomorphism \cite[Lemma 3.2.19]{AD}. Masuoka's Lemma
\ref{lema:iso-galois-extensions} implies another version of the
five-lemma: If $A$ and $A'$ are $H$-Galois over $B$, and $A'$ is
right $B$-faithfully flat, then $\theta$ is also an isomorphism.

\begin{cor}\label{cor:five}
Assume in \eqref{diag:five} that $\dim H$ is finite, $A'$ is
noetherian and $B$ is central in $A'$. Then $\theta$ is an
isomorphism.
\end{cor}

\pf As the rows are exact, the corresponding Galois maps $\beta$
and $\beta'$ are surjective; since $\dim H < \infty$, they are
bijective \cite[Thm. 1.7]{KT}. Thus $A$ and $A'$ are $H$-Galois
over $B$. Now $A'$ is $B$-faithfully flat by \cite[Thm.
3.3]{Sch1}. \epf

\section{Constructing quantum subgroups}
\label{section:q-subgroups}

In this section we construct quotients of the quantized coordinate
algebra $\Oe$. We do this in three steps.

\subsection{First step}\label{sec:fstep}
We construct in this subsection a quotient of $\Oe$ associated to a
Hopf subalgebra of $\qe$; it corresponds to a connected Lie subgroup
$L$ of $G$. Let $r: \qe^{*} \to H$ be a surjective Hopf algebra
morphism. Then we have an injective Hopf algebra map $\ ^{t}r: H^{*}
\to \qe$ and by Corollary \ref{cor:hopf-subalg-uepsilong}, the Hopf
algebra $H^{*}$ corresponds to a triple $(\Sigma, I_+, I_-)$. We
shall eventually show that this triple is part of a subgroup datum
as in Definition \ref{def:subgroupdatum}.

\subsubsection{The Hopf subalgebra $\Glil$ of $\Glie$}
\begin{definition}\label{def:gammaL}
For every triple $(\Sigma, I_+, I_-)$ define $\Gl$ to be the
subalgebra of $\Gli$ generated by the elements

\begin{align*}
K_{\alpha_{i}}^{-1} & & (1\leq i \leq n),\\
\left(\begin{array}{c}K_{\alpha_{i}}; 0\\
m\end{array}\right)&:= \prod_{s=1}^{m}\left(\frac{K_{\alpha_{i}}
q_{i}^{-s+1}-1}{q_{i}^{s}-1}\right) & (m\geq 1,\ 1\leq i \leq n),
\\
E_{j}^{(m)} & := \frac{E_{j}^{m}}{[m]_{q_{j}}!} &
 (m\geq 1,\ j\in I_+),\\
F_{k}^{(m)} &:= \frac{F_{k}^{m}}{[m]_{q_{k}}!} &
 (m\geq 1,\ k\in I_-),
\end{align*}

\noindent where $q_{i} = q^{d_{i}}$ for $1\leq i\leq n$. Note that
$\Gl$ does not depend on $\Sigma$.
\end{definition}

Choosing a reduced expression $s_{i_{1}}\cdots s_{i_{N}}$ of the
longest element of the Weyl group one can order totally the positive
part $\Phi_{+}$ of the root system $\Phi$ with $\beta_{1} =
\alpha_{i_{1}}$, $\beta_{2} = s_{i_{1}}\alpha_{i_{2}}, \ldots,\
\beta_{N} = s_{i_{1}}\cdots s_{i_{N-1}}\alpha_{i_{N}}$. Then using
the algebra automorphisms $T_{i}$ introduced by Lusztig \cite{L},
one may define corresponding root vectors $E_{\beta_{k}} =
T_{i_{1}}\cdots T_{i_{k-1}}E_{i_{k}}$ and $F_{\beta_{k}} =
T_{i_{1}}\cdots T_{i_{k-1}}F_{i_{k}}$. Consider now the
$R$-submodules of $\Gli$ given by
\begin{align*}
J_{\ell} = R \Big\{ \prod_{\beta\geq 0}F_{\beta}^{(n_{\beta})} \cdot
\prod_{i= 1}^{n}\left(\begin{array}{c}K_{\alpha_{i}}; 0\\
t_{i}\end{array}\right) K_{\alpha_{i}}^{\Ent(t_{i} / 2)} \cdot
\prod_{\alpha\geq 0}E_{\alpha}^{(m_{\alpha})} & : \\
\exists\ n_{\beta},t_{i}, m_{\alpha} \not\equiv 0 \mod (\ell) \Big\}\\
\\
\Ga_{\ell} = R \Big\{ \prod_{\beta\geq 0}F_{\beta}^{(n_{\beta})}
\cdot
\prod_{i= 1}^{n}\left(\begin{array}{c}K_{\alpha_{i}}; 0\\
t_{i}\end{array}\right) K_{\alpha_{i}}^{\Ent(t_{i} / 2)} \cdot
\prod_{\alpha\geq 0}E_{\alpha}^{(m_{\alpha})} & : \\
\forall\ n_{\beta},t_{i}, m_{\alpha} \equiv 0 \mod (\ell) \Big\}
\end{align*}

\noindent Then, by \cite[Thm. 6.3]{DL} there is a decomposition of
free $R$-modules $\Gli = J_{\ell}\ot \Ga_{\ell}$ and $ \Ga_{\ell}
/[p_{\ell}(q)\Ga_{\ell}] \simeq U(\g)_{\QQ(\epsilon)}$. Let
$Q_{I_{\pm}} = \bigoplus_{i\in I_{\pm}}^{} \ZZ \alpha_{i}$ and
define the following $R$-submodules of $\Gl$:
\begin{align*}
W_{\ell} =  R \Big\{& \prod_{\beta\geq 0}F_{\beta}^{(n_{\beta})}
\cdot
\prod_{i =1}^{n}\left(\begin{array}{c}K_{\alpha_{i}}; 0\\
t_{i}\end{array}\right) K_{\alpha_{i}}^{\Ent(t_{i} / 2)} \cdot
\prod_{\alpha\geq 0}E_{\alpha}^{(m_{\alpha})} : \\
& \exists\ n_{\beta},t_{i}, m_{\alpha} \not\equiv 0 \mod (\ell)
\mbox{ with }\beta \in Q_{I_{-}}, \alpha \in  Q_{I_{+}},
1\leq i \leq n \Big\}\\
\\
\Theta_{\ell} =  R \Big\{ & \prod_{\beta\geq 0}F_{\beta}^{(n_{\beta})}
\cdot
\prod_{i=1}^{n}\left(\begin{array}{c}K_{\alpha_{i}}; 0\\
t_{i}\end{array}\right) K_{\alpha_{i}}^{\Ent(t_{i} / 2)} \cdot
\prod_{\alpha\geq 0}E_{\alpha}^{(m_{\alpha})}  : \\
& \forall\ n_{\beta},t_{i}, m_{\alpha} \equiv 0 \mod (\ell)\mbox{
with }\beta \in  Q_{I_{-}}, \alpha \in  Q_{I_{+}},1\leq i \leq n
\Big\}
\end{align*}

\noindent Using the decomposition of $\Gli$ as free $R$-module we
get the following.

\begin{lema}\label{lema:decomp}
There is a decomposition of free $R$-modules $\Gl = W_{\ell}\ot
\Theta_{\ell}$. In particular, $\Gl$ is a direct summand of $\Gli$.
\end{lema}

\pf Clearly, $\Gl$ contains the free $R$-module $W_{\ell}\ot
\Theta_{\ell}$. Thus, it is enough to show that $\Gl \subseteq
W_{\ell}\ot \Theta_{\ell}$, but this follows directly from the fact
that $\Gl$ is generated as an algebra over $R$ by the elements in
Definition \ref{def:gammaL} and these generators satisfy the
relations given in \cite[Sec. 3.4]{DL}. \epf

Let $\Glil : = \Gl / [p_{\ell}(q)\Gl]$.
Then we have the following proposition.
\begin{prop}\begin{enumerate}
\item[$(a)$] $\Glil$ is a Hopf subalgebra of $\Glie$.
\item[$(b)$] $\Glie \simeq \Gli\ot_{R} R /[p_{\ell}(q)R]$
and $\Glil \simeq
\Gl\ot_{R} R /[p_{\ell}(q)R]$.
\end{enumerate}
\end{prop}

\pf We prove only $(a)$ since $(b)$ is straightforward. By
definition, the elements $E_{j}$ are
$(K_{\alpha_{j}},1)$-primitives, the $F_{k}$'s are
$(1,K_{\alpha_{k}}^{-1})$-primitives and the $K_{\alpha_{i}}$'s are
group-like. Moreover, the antipode is given by $\SS(K_{\alpha_{i}})
= K_{\alpha_{i}}^{-1}$, $\SS(E_{j}) = -K_{\alpha_{j}}^{-1}E_{j}$ and
$\SS(F_{k}) = -F_{k}K_{\alpha_{k}}$ with $1\leq i \leq n$, $j \in
I_{+}$ and $k\in I_{-}$. Hence, the subalgebra of $\Gl$ generated by
these elements is a Hopf subalgebra of $\Gli$ and $\Gl /
[p_{\ell}(q)\Gli\cap\Gl]$ is a Hopf subalgebra of $\Glie$. But by
Lemma \ref{lema:decomp}, we know that $\Gli = \Gl \oplus N$ for some
$R$-submodule $N$. Then $p_{\ell}(q)\Gli\cap\Gl = p_{\ell}(q)(\Gl
\oplus N)\cap\Gl = p_{\ell}(q)\Gl$, which implies that $\Glil  = \Gl
/ [p_{\ell}(q)\Gli\cap\Gl]$. \epf

\subsubsection{The regular Frobenius-Lusztig kernel $\ql$}\label{subsub:ul}
Let $\ql$ be the subalgebra of $\Glil$ generated by the elements
\begin{equation*}
\{K_{\alpha_{i}}, E_{j},F_{k}:\ 1\leq i \leq n, j\in I_{+}, k\in
I_{-}\}.
\end{equation*}

\begin{lema}\label{lema:uelsualg}
$\ql$ is a Hopf subalgebra of $\Glil$ such that $\Glil \cap \qe =
\ql$ and corresponds to the triple $(\mathbb{T},I_{+},I_{-})$, see
\eqref{eqn:torofinito}.
\end{lema}

\pf It is clear that $\ql$ is a Hopf subalgebra of $\Glil$. Since
the Frobenius-Lusztig kernel $\qe$ is the subalgebra of $\Glie$
generated by the elements $\{K_{\alpha_{i}}, E_{i},F_{i}:\ 1\leq
i\leq n\}$, we have that $\ql \subseteq \Glil \cap \qe$. But from
Lemma \ref{lema:decomp}, it follows that every element of $\Glil
\cap \qe$ must be contained in $\ql$. The last assertion follows
immediately from Corollary \ref{cor:hopf-subalg-uepsilong}. \epf

Recall that the quantum Frobenius map\label{qfrmap} $\Fr: \Glie \to
\Ege$ is defined on the generators of $\Glie$ by
\begin{align*}
\Fr(E_{i}^{(m)}) = \begin{cases}{\begin{matrix}e_{i}^{(m/\ell)} &
\mbox{if }
\ell\vert m\\
0& \mbox{otherwise,}\end{matrix}}\end{cases} & & \Fr(F_{i}^{(m)})
=
\begin{cases}\begin{matrix}f_{i}^{(m/\ell)} & \mbox{if } \ell\vert
m\\
0& \mbox{otherwise,}\end{matrix}\end{cases}\\
\Fr(\left(\begin{smallmatrix} K_{\alpha_{i}};0\\
m\end{smallmatrix}\right))
= \begin{cases}\begin{matrix} \left(\begin{smallmatrix} h_{i};0\\
m\end{smallmatrix}\right)& \mbox{if } \ell\vert
m\\
0& \mbox{otherwise,}\end{matrix}\end{cases} & &
\Fr(K_{\alpha_{i}}^{-1}) = 1,\mbox{ for all }1\leq i\leq n,
\end{align*}

\noindent and one has an exact sequence of Hopf algebras-- see
\cite{L}, \cite[Thm. 6.3]{DL}:
$$
1\to \qe \xrightarrow{} \Glie \xrightarrow{\Fr}
U(\lieg)_{\QQ(\epsilon)} \to 1.
  $$

\noindent If we define $\Ele := \Fr (\Glil)$, then it follows that
$\Ele$ is a subalgebra of $\Ege$ and the following diagram
commutes
\begin{equation}\label{diag:3}
\xymatrix{\qe\ar@{^{(}->}[r]^{}  &
\Glie\ar@{->>}[r]^{\Fr}& \Ege\\
\ql\ar@{^{(}->}[r]_{}\ar@{^{(}->}[u]_{} &\Glil
\ar@{->>}[r]^{\overline{\Fr}}\ar@{^{(}->}[u]_{} & \Ele,
\ar@{^{(}->}[u]_{}}
\end{equation}

\noindent where $\overline{\Fr}$ is the restriction of $\Fr$ to
$\Glil$.

\begin{rmk}\label{rmk:l}
$(a)$ Let $\liel$ be the set of primitive elements $P(\Ele)$ of
$\Ele$. Then $\liel$ is a Lie subalgebra of $\g$, which is in fact
regular in the sense of \cite{D1}: it is the Lie subalgebra
generated by the set $\{h_{i}, e_{j},f_{k}: 1\leq i \leq n, j\in
I_{+}, k\in I_{-}\}$. This agrees with Definition
\ref{def:subgroupdatum}.

\bigbreak $(b)$ $\Ker \overline{\Fr}$ is the two-sided ideal $\II$
of $\Glil$ generated by the set $$\Big\{E_{j}^{(m)}, F_{k}^{(m)},
\left(\begin{array}{c}K_{\alpha_{i}}; 0\\
m\end{array}\right), K_{\alpha_{i}} -1: 1\leq i \leq n, j\in
I_{+}, k\in I_{-}, m\geq 0, \ell \nmid m\Big\},$$

\noindent and coincides with $W_{\ell}$. Indeed, by \cite[Thm.
6.3]{DL} we know that $\Ker \Fr = J_{\ell}$ and coincides with the
two-sided ideal generated by
$$\Big\{E_{i}^{(m)}, F_{i}^{(m)},
\left(\begin{array}{c}K_{\alpha_{i}}; 0\\
m\end{array}\right), K_{\alpha_{i}} -1: 1\leq i\leq n, m\geq 0, \ell
\nmid m\Big\}.$$

\noindent But by Lemma \ref{lema:decomp}, $\Ker \overline{\Fr} =
\Ker \Fr \cap \Glil = J_{\ell} \cap \Glil = W_{\ell}$ and the last
one coincides with the ideal $\II$.

\bigbreak $(c)$ Since by \cite[Thm. 6.3]{DL}, the morphism $
\Ga_{\ell} /[p_{\ell}(q)\Ga_{\ell}] \to \Ege$ induced by the quantum
Frobenius map is bijective and by definition $\Theta_{\ell}
\subseteq \Ga_{\ell}$ and $\Ele = \Fr (\Ege)$, it follows by Lemma
\ref{lema:decomp} that $\Theta_{\ell}\cap p_{\ell}(q)\Ga_{\ell} =
p_{\ell}(q)\Theta_{\ell}$ and the morphism $ \Theta_{\ell}
/[p_{\ell}(q)\Theta_{\ell}] \to \Ele$ is also bijective.
\end{rmk}

The following proposition gives some properties of $\ql$.

\begin{prop}\label{prop:sucUl}
$(a)$ The following sequence of Hopf algebras is exact
\begin{equation}
1\to \ql \xrightarrow{j} \Glil \xrightarrow{\overline{\Fr}} \Ele \to
1.
\end{equation}
$(b)$ There is a surjective algebra map $\psi: \Glil \to \ql$ such
that $\psi|_{\ql} = \id$.
\end{prop}

\pf $(a)$ We need only to prove that $\Ker \overline{\Fr} =
\ql^{+}\Glil$ and $\ ^{\co \overline{\Fr}}\Glil = \ql$. The first
equality follows directly from Remark \ref{rmk:l} $(b)$, since the
two-sided ideal generated by $\ql^{+}$ coincides with $\II$. The
second equality follows from Lemma \ref{lema:uelsualg}, because $\
^{\co \Fr}\Glie = \qe$ by \cite[Lemma 3.4.1]{A} and $\ql = \qe \cap
\Glil = \ ^{\co \overline{\Fr}}\Glie \cap \Glil =\ ^{\co
\overline{\Fr}}\Glil$.

\bigbreak $(b)$ By Lemma \ref{lema:quotient}, there exists a
surjective algebra map $\varphi: \Glie \to \qe$ such that
$\varphi|_{\qe} = \id$. If we define $\psi: = \varphi|_{\ql}:
\Glil \to \qe$, then $\Img \psi \subseteq \ql$ and $\varphi|_{\ql}
= \id$, from which follows that $\Img \psi = \ql$. \epf

\subsubsection{The quantized coordinate algebra $\Oel$}
The inclusion $\Glil \hookrightarrow \Glie$ determines by duality a
Hopf algebra map $\Res: \Glie^{\circ} \to \Glil^{\circ}$. Since by
Proposition \ref{perfectpairing}, we have that $\Oeq \subseteq
\Glie^{\circ}$, we may define
$$\Oeql:=
\Res (\Oeq).$$

\noindent Moreover, as $\OG \subseteq \Oeq$, $\Res (\OG)$ is a
central Hopf subalgebra of $\Oeql$ and whence there exists an
algebraic subgroup $L$ of $G$ such that $\Res (\OG) = \OL$. Next
we show that $L$ is connected and the corresponding Lie subalgebra
of $\g$ is no other than the Lie algebra $\lgot$ discussed in
Remark \ref{rmk:l} (a).

\smallbreak Recall that a Lie subalgebra $\liek \subseteq \g$ is
called \emph{algebraic} if there exists an algebraic subgroup $K
\subseteq G$ such that $\liek = \Lie (K)$. We say that $\liek^{+}$
is the \emph{algebraic hull} of $\liek$ if $\liek^{+}$ is an
algebraic subalgebra of $\g$ such that $\liek \subseteq \liek^{+}$
and if $\agot$ is an algebraic subalgebra of $\lieg$ that contains
$\liek$, then $\liek^{+}\subseteq \agot$.

\begin{prop}
The algebraic group $L$ is connected and $\Lie(L) = \lgot$.
\end{prop}

\pf Since $\OG \subseteq \Ege^{\circ}$, dualizing diagram
\eqref{diag:3} we have $\OL $ $= \Res (\OG) \subseteq \Ele^{\circ}$.
But by \cite[XVI.3]{Ho}, $\Ele^{\circ}$ and consequently $\OL$ are
integral domains, implying that $L$ is irreducible and therefore
connected.

\smallbreak To show $\Lie(L) = \lgot$, we prove that $\Lie(L)$ is
the algebraic hull of $\lgot$ and $\lgot$ is an algebraic Lie
algebra. Since $\Ker \Res|_{\Oeq} = \{f\in \Oeq:\ f|_{\Glil} =0\}$
and the inclusion of $\OG$ in $\Oeq$ is given by the transpose of
the quantum Frobenius map $\Fr$ (see page \pageref{qfrmap}), it
follows that $\OL \simeq \OG / J$, where
\begin{align*}
J & = \{f\in \OG:\ \langle f,\Fr(x)\rangle = 0, \forall\ x\in \Glil\}\\
& = \{f\in \OG:\  \langle f,x\rangle = 0, \forall\ x\in \Ele\}.
\end{align*}

\noindent In particular,  $0 = \langle f,x\rangle = x(f)$ for all
$x\in \Ele$. Since by \cite[Lemma 6.9]{FR}, $\Lie(L) = \{\tau \in
\g:\ \tau(f) = 0,\ \forall f\in J\}$, it is clear that $\lgot
\subseteq \Lie(L)$. Now let $K \subseteq G$ such that $\lgot
\subseteq \Lie (K) =: \liek$ and denote by $\II$ the ideal of $K$;
then $\liek = \{\tau \in \g: \tau(\II) = 0\}$. As $\lgot \subseteq
\liek$, $\tau (\II) = 0$ for all $\tau \in \lgot$. Since the pairing
$\langle,\rangle$ is multiplicative, we have that $\II \subseteq J$
and whence $L \subseteq K$. Thus $\Lie (L) \subseteq \liek$ for all
algebraic Lie subalgebra $\liek$ such that $\lgot \subseteq \liek$,
implying that $\Lie(L) = \lgot^{+}$.

\smallbreak Now we show that $\lgot$ is algebraic, implying that
$\lgot = \lgot^{+} = \Lie(L)$. Consider $\g$ as a $G$-module with
the adjoint action and define $G_{\lgot} = \{x\in G:\ x\cdot \lgot =
\lgot\}$ and $\lieg_{\lgot} = \{\tau \in \g:\ [\tau, \lgot]
\subseteq \lgot\}$. Then by \cite[Ex. 8.4.7]{FR}, $\Lie(G_{\lgot}) =
\g_{\lgot}$. Thus, it is enough to show that $\lgot$ equals its
normalizer in $\g$.

\smallbreak By construction, we know that $\lgot = \lgot_{+} \oplus
\lieh \oplus \lgot_{-}$, where $\lieh$ is the Cartan subalgebra of
$\g$ and $\lgot_{\pm} = \bigoplus_{\alpha \in \Psi_{\pm}}
\g_{\alpha}$, with $\Psi_{\pm} = \{ \alpha \in \Phi:\ \sop(\alpha)
\subseteq I_{\pm}\}$. Let $x \in \g_{\lgot}$, then we may write $x=
\sum_{\alpha \in \Phi} c_{\alpha} x_{\alpha} + x_{0}$ with $x_{0}
\in \lieh$. Thus, for all $H \in \lieh$ we have that $[H, x] =
\sum_{\alpha \in \Phi} c_{\alpha} \alpha(H) x_{\alpha} \in \lgot$.
This implies that for all $H\in \lieh$, $c_{\alpha} \alpha(H) = 0$
for all $\alpha \notin \Psi = \Psi_{+} \cup \Psi_{-}$. Hence
$c_{\alpha} = 0$ for all $\alpha \notin \Psi$ and $x \in \lgot$.
\epf

\smallbreak Since $\OL$ is a central Hopf subalgebra of $\Oeql$, the
quotient
$$\overline{\Oel}_{\QQ(\epsilon)}:=  \Oeql /
[\OL^{+}\Oeql]$$ \noindent is a Hopf algebra which is
finite-dimensional. The following proposition shows that, as
expected, this algebra is isomorphic to $\ql^{*}$, see
\ref{subsub:ul}.

\begin{prop}\label{prop:sucOL}
$(a)$ The following sequence of Hopf algebras is exact
\begin{equation}\label{seq:obarra}
1\to \OL \xrightarrow{\iota_{L}} \Oeql \xrightarrow{\pi_{L}}
\overline{\Oel}_{\QQ(\epsilon)} \to 1.
\end{equation}
$(b)$ There exists a surjective Hopf algebra map $P: \qe^{*} \to
\overline{\Oel}_{\QQ(\epsilon)}$ making the following diagram
commutative:
\begin{equation}\label{diag:4}
\xymatrix{1 \ar[r]^{} & \OG\ar[r]^{\iota}\ar@{->>}[d]^{\res} &
\Oeq\ar[r]^{\pi}\ar@{->>}[d]^{\Res}& \qe^{*}\ar@{->>}[d]^{P}
\ar[r]^{} & 1\\
1 \ar[r]^{} & \OL\ar[r]^{\iota_{L}}  & \Oeql\ar[r]^{\pi_{L}}&
\overline{\Oel}_{\QQ(\epsilon)} \ar[r]^{} & 1.}
\end{equation}
\item[$(c)$] $\overline{\Oel}_{\QQ(\epsilon)} \simeq \ql^{*}$ as
Hopf algebras.
\end{prop}

\pf $(a)$ We need only to show that $\OL =\ ^{\co \pi_{L}}\Oeql$.
The algebra $\Oeq$ is noetherian, by Theorem \ref{OcentralenOe}
$(b)$. Therefore $\Oeql$ is also noetherian, since it is a quotient
of $\Oeq$. Then by \cite[Thm. 3.3]{Sch1}, $\Oeql$ is faithfully flat
over $\OL$ and by \cite[Prop. 3.4.3]{Mo} it follows that $\OL =\
^{\co \pi_{L}}\Oeql = \Oeql^{\co \pi_{L}}$.

\smallbreak $(b)$ Since the sequence \eqref{extOOeqe} is exact, we
have $\Ker \pi = \OG^{+}\Oeq$ and $\qe^{*} \simeq \Oeq /
[\OG^{+}\Oeq]$. But then, $\pi_{L}\Res (\Ker \pi) =
\pi_{L}(\OL^{+}\Oeql) =0$ and hence there exists a Hopf algebra map
$P : \qe^{*} \to \overline{\Oel}_{\QQ(\epsilon)}$ which makes the
diagram \eqref{diag:4} commutative.

\smallbreak $(c)$ Dualizing diagram \eqref{diag:3} we obtain a
commutative diagram
\begin{equation}\label{diag:5}
\xymatrix{\Ege^{\circ}\ar[d]_{}\ar@{^{(}->}[r]^{^{t}\Fr}  &
\Glie^{\circ}\ar[r]^{F}\ar[d]^{\Res} & \qe^{*}\ar@{->>}[d]^{p} \\
\Ele^{\circ}\ar@{^{(}->}[r]_{^{t}\overline{\Fr}}&\Glil^{\circ}
\ar[r]_{f} & \ql^{*}. }
\end{equation}

\noindent Since $\Oeql = \Res (\Oeq)$, $\OL = \Res (\OG)$ and
$\OG$ $ \simeq \Ege^{\circ}$, because $\g$ is simple, it follows
that $\OL \subseteq\ ^{t}\overline{\Fr}(\Ele^{\circ})$. In
particular, $\OL^{+} \subseteq \Ker f$. Moreover, since $F(\Oeq) =
\pi (\Oeq)$ $= \qe^{*}$ we have that $\ql^{*}= f\Res (\Oeq) =
f(\Oeql)$. Hence, there exists a surjective Hopf algebra map
$\beta :\overline{\Oel}_{\QQ(\epsilon)} \to \ql^{*}$; and $\dim
\overline{\Oel}_{\QQ(\epsilon)} \geq \dim \ql^{*}$.

\smallbreak We show next that there exists a surjective morphism
$\ql^{*} \to \overline{\Oel}_{\QQ(\epsilon)}$ implying that $\beta$
is an isomorphism. Consider the map $p: \qe^{*} \to \ql^{*}$ as in
\eqref{diag:5} and let $a \in \Ker p$. Since $\qe$ is
finite-dimensional, the coordinate functions of the regular
representation of $\qe$ span linearly $\qe^{*}$ and we may assume
that $a$ is a coordinate function of a finite-dimensional
representation $M$ of $\qe$. As $p$ is just the map given by the
restriction, we have that $a$ must be trivial on every basis of
$\ql$, in particular the following:
\begin{align*}
\Big\{ \prod_{\beta\geq 0}F_{\beta}^{n_{\beta}} \cdot
\prod_{i=1}^{n} K_{\alpha_{i}}^{t_{i}} \cdot \prod_{\alpha\geq
0}E_{\alpha}^{m_{\alpha}} & : 0\leq n_{\beta},t_{i}, m_{\alpha} <
\ell,
\\
& \beta \in Q_{I_{-}}, 1\leq i\leq n, \alpha \in Q_{I_{+}} \Big\}.
\end{align*}

\noindent On the other hand, we know by Lemma \ref{lema:quotient}
that there exists a surjective algebra map $\varphi: \Glie \to
\qe$ such that $\varphi|_{\qe} = \id$. Hence, the $\qe$-module $M$
admits a $\Glie$-module structure via $\varphi$. Since $M$ is
finite-dimensional and $K_{\alpha_{i}}^{\ell}$ acts as the
identity for every $1\leq i\leq n$, it follows that each operator
$K_{\alpha_{i}}$ is diagonalizable with eigenvalues
$\epsilon_{i}^{m}$ for some $m\in \NN$. This implies by definition
that the coordinate function $\varphi^{*}(a)$ of the
$\Glie$-module $M$ must be contained in $\Oeq$. Thus, using the
definition of $\varphi$ we have that $\Res \varphi^{*}(a)$ must
annihilate the set
\begin{align*}W_{\ell} =  \QQ(\epsilon)
\Big\{& \prod_{\beta\geq 0}F_{\beta}^{(n_{\beta})} \cdot
\prod_{i =1}^{n}\left(\begin{array}{c}K_{\alpha_{i}}; 0\\
t_{i}\end{array}\right) K_{\alpha_{i}}^{\Ent(t_{i} / 2)} \cdot
\prod_{\alpha\geq 0}E_{\alpha}^{(m_{\alpha})} : \\
& \exists\ n_{\beta},t_{i}, m_{\alpha} \not\equiv 0 \mod (\ell)
\mbox{ with }\beta \in Q_{I_{-}}, 1\leq i\leq n, \alpha \in
Q_{I_{+}} \Big\}.
\end{align*}

\noindent Since by Lemma \ref{lema:decomp}, $\Gl = W_{\ell}\ot
\Theta_{\ell}$ as free $R$-modules and by Remark \ref{rmk:l}, $\Ker
\overline{\Fr} = W_{\ell}$ and the map $\Theta_{l}/
[p_{\ell}(q)\Theta_{\ell}]\to \Ele$ induced by the restriction of
the quantum Frobenius map $\overline{\Fr}$ is bijective. Then there
exists $b\in \Ele^{\circ}$ such that $^{t}\overline{\Fr}(b) = \Res
(\varphi^{*}(a))$. Hence,
$$P(a) = P(\pi(\varphi^{*}(a))) =
\pi_{L}(\Res(\varphi^{*}(a))) =\pi_{L}(^{t}\overline{\Fr}(b)) =
\eps(b) = \eps(a) = 0,$$

\noindent and $a \in \Ker P$. Thus $\Ker p \subseteq \Ker P$ and
there exists a surjective map $\ql^{*} \to
\overline{\Oel}_{\QQ(\epsilon)}$. \epf

\begin{obs}\label{rmk:diagcomm}
By the Proposition above, we have the following commutative
diagram of exact sequences of Hopf algebras

\begin{equation}\label{diag:6}
\xymatrix{1\ar[r]^{}& \OG\ar[r]^{\iota}\ar@{->>}[d]^{\res} &
\Oeq\ar[r]^{\pi}\ar@{->>}[d]^{\Res}& \qe^{*}\ar@{->>}[d]^{p}\ar[r]^{}& 1\\
1\ar[r]^{} & \OL\ar[r]^{\iota_{L}}  & \Oeql\ar[r]^{\pi_{L}}& \ql^{*}
\ar[r]^{} & 1}
\end{equation}
\end{obs}

\subsection{Second Step}
We consider now the complex form of the algebras defined above.
Denote the $\CC$-form of the Frobenius-Lusztig kernels just by $\qe$
and $\ql$.

The following proposition tell us how to construct Hopf algebras
from a central exact sequence and a surjective Hopf algebra map.
We perform it in a general setting and then we apply it to our
situation. The characterization of these algebras as pushouts will
be crucial.

\begin{prop}\label{cociente1}
Let $A$ and $K$ be Hopf algebras, $B$ a central Hopf subalgebra of
$A$ such that $A$ is left or right faithfully flat over $B$ and
$p: B \to K$ a surjective Hopf algebra map. Then $H = A/AB^{+}$ is
a Hopf algebra and $A$ fits into the exact sequence $1\to B
\xrightarrow{\iota} A \xrightarrow{\pi} H \to 1$. If we set $\JJ =
\Ker p \subseteq B$, then $(\JJ) = A\JJ $ is a Hopf ideal of $A$
and $A / (\JJ)$ is the pushout given by the following diagram:

$$\xymatrix{B\,
\ar@{^{(}->}[r]^{\iota} \ar@{->>}[d]_{p} & A \ar@{->>}[d]^{q}\\
K\ar@{^{(}->}[r]_(.4){j} & A / (\JJ) .}$$

\noindent Moreover, $K$ can be identified with a central Hopf
subalgebra of $A / (\JJ)$ and $A / (\JJ)$ fits into the exact
sequence
\begin{equation}\label{eqn:exactseq-pushout}
1\to K \to A/(\JJ) \to H \to 1.
\end{equation}
\end{prop}

\pf The first assertion follows directly from \cite[Prop.
3.4.3]{Mo}. Since $B$ is central in $A$, $(\JJ)$ is a two-sided
ideal of $A$. Moreover, from the fact that $\eps$ and $\com$ are
algebra maps and $\SS (\JJ) \subseteq \JJ$, it follows that $(\JJ)$
is indeed a Hopf ideal. Identify $K$ with $B/\JJ$. Then the map $j:
K \to A/(\JJ)$ given by $j(b + \JJ) = \iota(b) + (\JJ)$ defines a
morphism of Hopf algebras because $\iota$ is a Hopf algebra map.
Since $A$ is faithfully flat over $B$, by \cite[Cor. 1.8]{Sch2}, $B$
is a direct summand in $A$ as a $B$-module, say $A = B \oplus M$.
Then $(\JJ) \cap B = \JJ A \cap B = (\JJ B \oplus \JJ M)\cap B$ $=
(\JJ \oplus \JJ M)\cap B = \JJ$. Thus, if $j(b + \JJ)= 0$ then
$\iota (b) \in (\JJ)$ and this implies that $b \in (\JJ)\cap B =
\JJ$ by the equality above. Hence, $j$ is injective.

 Let us see now that $A/(\JJ)$ is a pushout: let $C$ be a
Hopf algebra and suppose that there exist Hopf algebra maps
$\varphi_{1} : K \to C$ and $\varphi_{2} : A \to C$ such that
$\varphi_{1} p = \varphi_{2}\iota $. We have to show that there
exists a unique Hopf algebra map $\phi: A / (\JJ) \to C$ such that
$\phi q = \varphi_{2}$ and $\phi j = \varphi_{1}$.

$$\xymatrix{B
\ar[r]^{\iota} \ar[d]_{p} &
A \ar[d]^{q} \ar@/^/[ddr]^{\varphi_{2}}\\
K\ar[r]_{j} \ar@/_/[drr]_{\varphi_{1}}& A / (\JJ)
\ar@{-->}[dr]|-{\exists ! \phi}\\
& & C}$$

\smallbreak Since $\varphi_{2} ((\JJ)) = \varphi_{2} (A\JJ) =
\varphi_{2}(A) \varphi_{2} (\iota (\JJ)) = \varphi_{2}(A)
\varphi_{1} (p(\JJ))= 0$, there exists a unique Hopf algebra map
$\phi: A / (\JJ) \to C$ such that $\phi q = \varphi_{2}$. Moreover,
let $x \in K$ and $b \in B$ such that $p(b) = x$. Then $\phi  j (x)
= \phi j p (b)$ $= \phi q \iota (b) $ $= \varphi_{2} \iota (b) $ $=
\varphi_{1} p (b) = \varphi_{1} (x)$, from which follows that $\phi
j = \varphi_{1}$.

\smallbreak Denote also by $K$ the image of $K$ under $j$. To see
that $K$ is central in $A/(\JJ)$ we have to verify that $
j(c)\bar{a}  = \bar{a}j(c)$ for all $\bar{a} \in A/(\JJ)$, $c \in
K$. Since $p$ is surjective, for all $c \in K$ there exists $b \in
B$ such that $p(b) = c$ and since $q$ is an algebra map, it follows
that $\bar{a}j(c) = q(a)j(p(b))= $ $q(a)q(\iota(b)) = q(a\iota(b))$
$= q(\iota(b)a) = q(\iota(b))q(a) = j(c)\bar{a}$, because $B$ is
central in $A$. In particular, the quotient $\widetilde{H}$ $=
[A/(\JJ)]/[K^{+}(A/(\JJ))]$ is a Hopf algebra. To see that $A/(\JJ)$
is a central extension of $K$ by $\widetilde{H}$, by \cite[Prop.
3.4.3]{Mo} it is enough to show that $A/(\JJ)$ is flat over $K$ and
$K$ is a direct summand of $A/(\JJ)$ as $K$-modules, since by
\cite[Cor. 1.8]{Sch2} this implies that $A/(\JJ)$ is faithfully flat
over $K$.

\smallbreak First we show that $A/(\JJ)$ is flat over $K$. Let
$M_{1}$ and $M_{2}$ be two right $K$-modules and let $f: M_{1} \to
M_{2}$ be an injective homomorphism. In particular, they admit a
$B$-module structure via the map $p: B \to K$, which we denote by
$\overline{M}_{i}$ for $i = 1,2$; thus $f$ is an injective
homomorphism of $B$-modules. Since $A$ is faithfully flat over $B$,
the homomorphism of $A$-modules $f\otimes \id :
\overline{M}_{1}\otimes_{B} A \to \overline{M}_{2}\otimes_{B} A$ is
also injective. As $\JJ$ is central in $A$, we have for $i=1,2$ that
$(\overline{M}_{i}\otimes_{B} A) (\JJ) = 0$. Then the $A$-modules
are also $A/(\JJ)$-modules and $\overline{M}_{i}\otimes_{B} A \simeq
M_{i}\otimes_{K} A/(\JJ)$ as $A/(\JJ)$-modules by the construction
of $\overline{M}_{i}$. Hence the homomorphism of $A/(\JJ)$-modules
$f\otimes \id: M_{1}\otimes_{K} A/(\JJ) \to M_{2}\otimes_{K}
A/(\JJ)$ is injective and $A/(\JJ)$ is flat over $K$.

\smallbreak As $A = B \oplus M$ as $B$-modules, we have that $(\JJ)
= A\JJ = \JJ \oplus M\JJ$, where $M\JJ$ is a $B$-submodule of $M$
and $\JJ = B\cap (\JJ \oplus M\JJ)$. Hence $A/(\JJ) = (B\oplus M)
/(\JJ \oplus M\JJ) = K \oplus (M/M\JJ)$ as $K$-modules, which
implies that $K$ is a direct summand of $A/(\JJ)$.

\smallbreak In conclusion, $A/(\JJ)$ fits into an exact sequence of
Hopf algebras $$ 1\to K \xrightarrow{j} A/(\JJ) \xrightarrow{r}
\widetilde{H} \to 1. $$

\noindent Since the map $\Psi: K^{+}(A/(\JJ)) \to (B^{+}A) /(\JJ)$
defined by $\Psi(\overline{b}\overline{a}) = \overline{ba}$ is a
$k$-linear isomorphism, it follows that $\widetilde{H} = (A /(\JJ))
/ [K^{+}(A/(\JJ))]$ $ \simeq (A/(\JJ)) / [(B^{+}A) /(\JJ)]$ $\simeq
A /B^{+}A = H$ and therefore $A/(\JJ)$ fits into an exact sequence
\eqref{eqn:exactseq-pushout}. \epf

\smallbreak Let $\Ga$ be an algebraic group and let $\sigma: \Ga \to
G$ an injective homomorphism of algebraic groups such that $\sigma
(\Ga) \subseteq L$. Then we have a surjective Hopf algebra map $\
^{t}\sigma: \OO(L) \to \OO(\Ga)$. Applying the pushout construction
given in Proposition \ref{cociente1}, we obtain a Hopf algebra $A_{
\lgot, \sigma}$ which is part of an exact sequence of Hopf algebras
and fits into the following commutative diagram

\begin{equation}\label{diag:pushoutol}
\xymatrix{1 \ar[r]^{} & \Oc(G) \ar[r]^{\iota} \ar[d]_{\res} & \Oe
\ar[r]^{\pi}
\ar[d]^{\Res} & \qe^{*} \ar[r]^{}\ar[d]^{p} & 1\\
1 \ar[r]^{} &   \Oc(L)\ar[d]_{^{t}\sigma}  \ar[r]^{\iota_L} & \Ol
\ar[r]^{\pi_L}\ar[d]^{\nu} &
\ql^* \ar[r]^{}\ar@{=}[d] & 1\\
1 \ar[r]^{} &   \OO(\Ga)  \ar[r]^{j} & A_{\lgot,
\sigma}\ar[r]^{\bar{\pi}} & \ql^* \ar[r]^{} & 1.}
\end{equation}

\begin{obs}\label{rmk:kt}
Let $1\to K \to A \to H\to 1$ be an exact sequence of Hopf algebras.
If $\beta : A\otimes_{K} A\to A\otimes H$, $\beta(x, y) =
xy_{(0)}\ot y_{(1)}$ denotes the Galois map, then $\beta$ is
surjective, since $H \simeq A /K^{+}A$. If moreover $H$ is
finite-dimensional, $A$ is a finitely generated projective
$K$-module, by \cite[Thm. 1.7]{KT}. In particular, if $\dim K$ is
finite, then $\dim A = \dim K \dim H$ is also finite. In our case,
if $\Ga$ is finite we obtain that $\dim A_{\lgot,\sigma}= \vert
\Ga\vert \dim \ql$.
\end{obs}

\subsection{Third Step}
In this subsection we make the third and last step of the
construction. It consists essentially on taking a quotient by a Hopf
ideal generated by differences of central group-like elements of
$A_{\lgot, \sigma}$. The crucial point here is the description of
$H$ as a quotient of $\ql^{*}$ and the existence of a coalgebra
morphism $\psi^{*} : \ql^{*} \to \Ol$.

\smallbreak Recall that from the beginning of this section we fixed
a surjective Hopf algebra map $r: \qe^{*} \to H$ and $H^{*}$ is
determined by the triple $(\Sigma, I_{+}, I_{-})$. Since the Hopf
subalgebra $\ql$ is determined by the triple $(\mathbb{T},I_{+},
I_{-})$ with $\mathbb{T} \supseteq \Sigma$, we have that $H^{*}
\subseteq \ql \subseteq \qe$. Denote by $v : \ql^{*} \to H$ the
surjective Hopf algebra map induced by this inclusion. Then $H$ is a
quotient of $\ql^{*}$ which fits into the following commutative
diagram
$$\xymatrix{\qe^{*} \ar@{->>}[r]^{p}\ar@{->>}[dr]_{r} &
\ql^{*}\ar@{->>}[d]^{v}\\
& H.  }$$

\begin{obs}\label{rmk:subgruposigma}
Let $I = I_{+} \cup -I_{-}$, $I^{c} = \Pi - I$  and $\mathbb{T}_{I}=
\{K_{\alpha_{i}}:\ i\in I\}$. Let $s=\vert I^{c} \vert$. By
Corollary \ref{cor:hopf-subalg-uepsilong}, we know that
 $\mathbb{T}_{I} \subseteq \Sigma \subseteq \mathbb{T} =
 \mathbb{T}_{I}\times \mathbb{T}_{I^{c}}$. If we set $\Omega
 = \Sigma \cap \mathbb{T}_{I^{c}}$, it follows clearly that
 $\Sigma \simeq \mathbb{T}_{I} \times \Omega$.

\smallbreak Thus, giving a subgroup $\Sigma$ such that
$\mathbb{T}_{I} \subseteq \Sigma \subseteq \mathbb{T}$ is the same
as giving a subgroup $\Omega \subseteq \mathbb{T}_{I^{c}}$, and this
is the same as giving a subgroup $N \subseteq
\widehat{\mathbb{T}_{I^{c}}}$. Namely, $N$ is the kernel of the
group homomorphism $\rho: \widehat{\mathbb{T}_{I^{c}}} \to
\widehat{\Omega}$ induced by the inclusion. In particular, we have
that $ \vert \Sigma \vert = \vert\mathbb{T}_{I}\vert
\vert\Omega\vert = \ell^{n-s}\vert \Omega \vert =
\frac{\ell^{n}}{\vert N\vert}$.
\end{obs}

\begin{definition}\label{def:D}
For all $1\leq i \leq n$ such that $\alpha_{i} \notin I_{+}$ or
$\alpha_{i}\notin I_{-}$ we define $D_{i} \in G(\ql^{*}) =
\Alg(\ql, \CC)$ on the generators of $\ql$ by
\begin{align*}
D_{i} (E_{j}) & = 0 \qquad \forall\ j:\ \alpha_{j} \in I_{+}, &
D_{i} (F_{k}) & = 0 \qquad \forall\ k:\ \alpha_{k} \in I_{-},\\
D_{i} (K_{\alpha_{t}}) & = 1 \qquad \forall\ t \neq i,\ 1\leq t\leq
n, & D_{i} (K_{\alpha_{i}}) & = \epsilon_{i},
\end{align*}
\noindent where $\epsilon_{i}$ is a primitive $\ell$-th root of 1.
If $\alpha_{i}\notin I_{+}$ or $\alpha_{i}\notin I_{-}$, then
$E_{i}$ or $F_{i}$ is not a generator of $\ql$, respectively.
Hence, $D_{i}$ is a well-defined algebra map, since it verifies
all the defining relations of $\Glie$ \cite[Sec. 3.4]{DL}, see
\cite[5.2.12]{G} for details.

\smallbreak Let $I^{c} = \{\alpha_{i_{1}},\ldots, \alpha_{i_{s}}\}$
and let $N \subseteq \widehat{\mathbb{T}_{I^{c}}}$, correspond to
$\Sigma$ as in Remark \ref{rmk:subgruposigma}. We define for all $z
= (z_{1},\ldots, z_{s}) \in \widehat{\mathbb{T}_{I^{c}}}$ the following
group-like element
$$D^{z}:= D_{i_{1}}^{z_{1}}\cdots D_{i_{s}}^{z_{s}}.$$
\end{definition}

Recall that $(M)$ denotes the two-sided ideal generated by a
subset $M$ of an algebra $R$.

\begin{lema}\label{lema:elemcentrales}
\begin{enumerate}
\item[$(a)$] If $\alpha_{i} \in I^{c}$ then $D_{i}$ is central in $\ql^{*}$.
In particular $D^{z}$ is central for all $z\in
\widehat{\mathbb{T}_{I^{c}}}$.
\item[$(b)$] $H \simeq \ql^{*} / (D^{z} - 1\vert z\in N)$.
\end{enumerate}
\end{lema}

\pf $(a)$ We have to show that $D_{i} f = f D_{i}$ for all $f \in
\ql^{*}$. First observe that $D_{i}$ coincide with the counit of
$\ql$ in all elements of the basis which do not contain some
positive power of $K_{\alpha_{i}}$. By Lemma \ref{lema:decomp} we
know that $\ql$ has a basis of the form
\begin{align*}
\Big\{ \prod_{\beta\geq 0}F_{\beta}^{n_{\beta}} \cdot \prod_{i
=1}^{n}K_{\alpha_{i}}^{t_{i}} \cdot \prod_{\alpha\geq
0}E_{\alpha}^{m_{\alpha}} : &\quad 0\leq n_{\beta}, t_{i},
m_{\alpha}<
\ell,\\
& \mbox{ with }\beta \in Q_{I_{-}}, \alpha \in Q_{I_{+}}, 1\leq i
\leq n \Big\}.
\end{align*}

\noindent Thus, using the defining relations of $\Glie$ \cite[Sec.
3.4]{DL}, we may assume that this basis is of the form
$K_{\alpha_{i}}^{t_{i}} M$ with $0\leq t_{i} < \ell$ and $M$ does
not contain any power of $K_{\alpha_{i}}$. Then for every element
of this basis we have
\begin{align*}
D_{i} f(K_{\alpha_{i}}^{t_{i}} M) & = D_{i} (K_{\alpha_{i}}^{t_{i}}
M_{(1)})f(K_{\alpha_{i}}^{t_{i}} M_{(2)})
=D_{i} (K_{\alpha_{i}}^{t_{i}}) D_{i} (M_{(1)})f(K_{\alpha_{i}}^{t_{i}} M_{(2)})\\
& = \epsilon_{i}^{t_{i}} \eps (M_{(1)})f(K_{\alpha_{i}}^{t_{i}}
M_{(2)}) =
\epsilon_{i}^{t_{i}} f(K_{\alpha_{i}}^{t_{i}} M) \\
& = f D_{i} (K_{\alpha_{i}}^{t_{i}} M),
\end{align*}

$(b)$ By $(a)$ we know that $D^{z}$ is a central group-like
element of $\ql^{*}$ for all $z \in N$. Hence the quotient
$\ql^{*} / (D^{z} - 1\vert z\in N)$ is a Hopf algebra.

\smallbreak On the other hand, following Corollary
\ref{cor:hopf-subalg-uepsilong} we know that $H^{*}$ is determined
by the triple $(\Sigma, I_{+},I_{-})$ and consequently $H^{*}$ is
included in $\ql$. If we denote $v: \ql^{*} \to H$ the surjective
map induced by this inclusion, we have that $\Ker v = \{f\in
\ql^{*}:\ f(h) = 0,\ \forall\ h\in H^{*}\}$. But $D^{z} -1 \in
\Ker v$ for all $z\in N$, since $D^{z}(\omega) = \rho(z)(\omega) =
1$ for all $\omega \in \Omega$. Hence there exists a surjective
Hopf algebra map
$$\gamma: \ql^{*}/(D^{z}-1\vert\
z\in N) \twoheadrightarrow H.$$

\noindent Combining Corollary \ref{cor:hopf-subalg-uepsilong} with
the PBW-basis of $H$ and $\ql$ we have that
\begin{align*}
\dim H & = \ell^{\vert I_{+}\vert + \vert I_{-} \vert}\vert \Sigma
\vert = \ell^{\vert I_{+}\vert + \vert I_{-} \vert} \ell^{n-s} \vert
\Omega \vert = \ell^{\vert I_{+}\vert + \vert I_{-} \vert}
\ell^{n-s} \vert \widehat{\Omega} \vert  = \ell^{\vert I_{+}\vert +
\vert I_{-} \vert}
\frac{\ell^{n}}{\vert N \vert}\\
& = \dim (\ql^{*}/(D^{z}-1\vert\ z\in N)),
\end{align*}

\noindent which implies that $\gamma$ is an isomorphism.
\epf

\begin{obs}
The lemma above is very similar to a result used by E. M\"uller in
the case of type $A_{n}$ \cite[Sec. 4]{Mu2} for the classification
of the finite-dimensional quotients of $\Oesln$. The new point of
view here consists in regarding $H$ as a quotient of the dual of
$\ql$.
\end{obs}

Before going on with the construction we need the following
technical lemma. Let $\x = \{D^{z}\vert\ z\in
\widehat{\mathbb{T}_{I^{c}}}\}$ be the set of central group-like
elements of $\ql^{*}$ given by Lemma \ref{lema:elemcentrales}.

\begin{lema}\label{lema:grupoz}
There exists a subgroup $\z:=\{\partial^{z}\vert\ z \in
\widehat{\mathbb{T}_{I^{c}}}\}$ of $G(A_{\lgot, \sigma})$ isomorphic
to $\x$ consisting of central elements.
\end{lema}

\pf By Proposition \ref{prop:sucUl} $(b)$, we know that there exists
an algebra map $\psi: \Glil \to \ql$; it induces a coalgebra map
$\psi^{*}: \ql^{*} \to \Glil^{\circ}$ such that the following
diagram commutes
$$\xymatrix{\Glie^{\circ} \ar@{->>}[d]_{\Res} & \qe^{*} \ar[l]_{\varphi^{*}}
\ar@{->>}[d]^{p}\\
\Glil^{\circ} & \ql^{*}. \ar[l]^{\psi^{*}}}$$

\noindent Here, $\varphi^{*}$ is the coalgebra map induced by the
algebra map $\varphi: \Glie \to \ql$ given by Lemma
\ref{lema:quotient}, whose restriction to $\Glil$ defines $\psi$.
Furthermore, by  the proof of Proposition \ref{prop:sucUl} $(c)$,
$\Img \varphi^{*}\subseteq \Oe$; since $\Res(\Oe) = \Ol$, it follows
that $\Img \psi^{*}\subseteq \Ol$. Consequently, we obtain a group
of group-like elements $\y = \{d^{z} = \psi^{*}(D^{z})\vert\ z\in
\widehat{\mathbb{T}_{I^{c}}}\}$ in $\Ol$. Moreover, by Lemma
\ref{lema:decomp} and the definitions of $\psi$ and the elements
$D_{i}$, the elements of $\y$ are central.

\smallbreak Since the map $\nu: \Ol \to A_{\lgot, \sigma}$ given by
the pushout construction is surjective, the image of $\y$ defines a
group of central group-like elements in $A_{\lgot, \sigma}$:
$$\z = \{\partial^{z} = \nu(d^{z})\vert\ z\in
\widehat{\mathbb{T}_{I^{c}}}\}.$$

\noindent Besides, $\vert \z \vert = \vert \y \vert = \vert \x
\vert = \ell^{s}$. Indeed, $\bar{\pi}(\z) = \bar{\pi}\nu(\y) =
\pi_{L}(\y) = \pi_{L}\psi^{*}(\x)= \x$ since the diagram
\eqref{diag:pushoutol} is commutative and $\pi_{L}\psi^{*} = \id$.
Hence $\vert \bar{\pi}(\z) \vert = \vert\x \vert$, from which the
assertion follows. \epf

We are now ready for our first main result.

\begin{theorem}\label{teo:constrfinal}
Let $\D = (I_{+}, I_{-}, N, \Ga, \sigma, \delta)$ be a subgroup
datum. Then there exists a Hopf algebra $A_{\D}$ which is a
quotient of $\Oe$ and fits into the exact sequence
$$ 1\to \Oc(\Ga) \xrightarrow{\hat{\iota}} A_{\D}
\xrightarrow{\hat{\pi}} H \to 1. $$

\noindent Concretely, $A_{\D}$ is given by the quotient $A_{\lgot,
\sigma}/J_{\delta}$ where $J_{\delta}$ is the two-sided ideal
generated by the set $\{\partial^{z} - \delta(z)\vert z\in N\}$
and the following diagram of exact sequences of Hopf algebras is
commutative

\begin{equation}\label{diag:const}
\xymatrix{1 \ar[r]^{}
& \Oc(G) \ar[r]^{\iota} \ar[d]_{\res} & \Oe \ar[r]^{\pi}
\ar[d]^{\Res} & \qe^{*} \ar[r]^{}\ar[d]^{p} & 1\\
1 \ar[r]^{} &   \Oc(L)\ar[d]_{^{t}\sigma}  \ar[r]^{\iota_L} & \Ol
\ar[r]^{\pi_L}\ar[d]^{\nu} &
\ql^* \ar[r]^{}\ar@{=}[d] & 1\\
1 \ar[r]^{} &   \Oc(\Ga)  \ar[r]^{j}\ar@{=}[d] & A_{\lgot,
\sigma}\ar[r]^{\bar{\pi}}\ar[d]^{t} & \ql^* \ar[r]^{}\ar[d]^{v} & 1\\
1 \ar[r]^{} &   \Oc(\Ga)  \ar[r]^{\hat{\iota}} &
A_{\D}\ar[r]^{\hat{\pi}}& H\ar[r]^{}& 1.}\end{equation}
\end{theorem}

\pf By Remark \ref{rmk:subgruposigma}, $N$ determines a subgroup
$\Sigma$ of $\mathbb{T}$ and the triple $(\Sigma, I_{+}, I_{-})$
give rise to a surjective Hopf algebra map $r: \qe^{*} \to H$.
Since $\sigma: \Ga \to L \subseteq G$ is injective, by the first
two steps developed before one can construct a Hopf algebra
$A_{\lgot,\sigma}$ which is a quotient of $\Oe$ and an extension
of $\Oc(\Ga)$ by $\ql^{*}$, where $\ql$ is the Hopf subalgebra of
$\qe$ associated to the triple $(\mathbb{T},I_{+},I_{-})$.
Moreover, by Lemma \ref{lema:elemcentrales} $(b)$, $H$ is the
quotient of $\ql^{*}$ by the two-sided ideal $(D^{z}-1\vert\ z\in
N)$. If $\delta: N \to \widehat{\Ga}$ is a group map, then the
elements $\delta(z)$ are central group-like elements in $A_{\lgot,
\sigma}$ for all $z \in N$, and the two-sided ideal $J_{\delta}$
of $A_{\lgot, \sigma}$ generated by the set $\{\partial^{z} -
\delta(z)\vert z\in N\}$ is a Hopf ideal. Hence, by \cite[Prop.
3.4 (c)]{Mu2} the following sequence is exact
$$ 1\to \Oc(\Ga)/ \J \to A_{\lgot,
\sigma}/J_{\delta} \to \ql^{*} /\bar{\pi}(\J_{\delta})\to 1,$$

\noindent where $\J = J_{\delta} \cap \Oc(\Ga)$. Since
$\bar{\pi}(\partial^{z}) = D^{z}$ and $\bar{\pi}(\delta(z)) = 1$ for
all $z\in N$, we have that $\bar{\pi}(J_{\delta})$ is the two-sided
ideal of $\ql^{*}$ given by $(D^{z} -1\vert\ z\in N)$, which implies
by Lemma \ref{lema:elemcentrales} $(b)$ that $\ql^{*} /
\bar{\pi}(\J_{\delta}) = H$. Hence, if we denote $A_{\D} : =
A_{\lgot, \sigma}/J_{\delta}$, we can re-write the exact sequence of
above as
\begin{equation}\label{suc:casi} \xymatrix{1\to \Oc(\Ga)/ \J \to
A_{\D} \to H \to 1.}\end{equation}

To end the proof it is enough to see that $\J = J_{\delta} \cap
\Oc(\Ga)= 0$. Clearly, $J_{\delta}$ coincides with the two-sided
ideal $(\partial^{z}\delta(z)^{-1}-1\vert\ z\in N)$ of
$A_{\lgot,\sigma}$. Moreover,  $\Upsilon :=
\{\partial^{z}\delta(z)^{-1}\vert\ z\in N\}$ is a subgroup of
central group-like elements of $G(A_{\lgot,\sigma})$ and $J_{\delta}
= (g-1\vert\ g\in \Upsilon) = A_{\lgot,\sigma}\CC[\Upsilon]^{+}$.
Let $\partial N = \{\partial^{z} \vert\ z\in N\}$. Then clearly the
subalgebra $B := \Oc(\Ga)\CC[\partial N]$ is a central Hopf
subalgebra of $A_{\sigma}$ which contains $\CC[\Upsilon]$. Further,
$B \simeq \Oc(\tilde{\Ga})$ for some algebraic group $\tilde{\Ga}$
and one has the following exact sequence of Hopf algebras
$$ 1\to \Oc(\Ga) \to \Oc(\tilde{\Ga}) \to R \to 1,$$

\noindent where $R = \Oc(\tilde{\Ga}) /
\Oc(\tilde{\Ga})\Oc(\Ga)^{+}$. But $R \simeq
\bar{\pi}(\Oc(\tilde{\Ga})) = \CC[N]$, since
\begin{align*}
\bar{\pi}(\Oc(\tilde{\Ga})) & = [\Oc(\tilde{\Ga}) +
\Oc(\Ga)^{+}A_{\lgot,\sigma}]/[\Oc(\Ga)^{+}A_{\lgot,\sigma}] \simeq
\Oc(\tilde{\Ga}) / [\Oc(\tilde{\Ga}) \cap
(\Oc(\Ga)^{+}A_{\lgot,\sigma})]\\
& \simeq \Oc(\tilde{\Ga}) / \Oc(\tilde{\Ga})\Oc(\Ga)^{+}.
\end{align*}

\noindent The last isomorphism follows from the fact that
$\Oc(\tilde{\Ga}) \cap (\Oc(\Ga)^{+}A_{\lgot,\sigma}) =
\Oc(\tilde{\Ga})\Oc(\Ga)^{+}$. Indeed, since $\Oc(\tilde{\Ga})$ is
a central Hopf subalgebra of the noetherian algebra
$A_{\lgot,\sigma}$, by \cite[Thm. 3.3]{Sch2}, $\Oc(\tilde{\Ga})$
is a direct summand of $A_{\lgot,\sigma}$ as
$\Oc(\tilde{\Ga})$-module, say $A_{\lgot,\sigma} =
\Oc(\tilde{\Ga}) \oplus M$. Then $\Oc(\Ga)^{+}A_{\lgot,\sigma} =
\Oc(\Ga)^{+}\Oc(\tilde{\Ga}) \oplus \Oc(\Ga)^{+}M$ and the claim
follows since $\Oc(\tilde{\Ga})\cap \Oc(\Ga)^{+}M = 0$. Hence we
have an exact sequence
$$ 1\to \Oc(\Ga) \to \Oc(\tilde{\Ga}) \xrightarrow{\bar{\pi}}
\CC[N] \to 1,$$

\noindent which is cleft by the proof of Lemma \ref{lema:grupoz},
since $\bar{\pi}$ admits a coalgebra section. Moreover, this section
on $\CC[N]$ is by definition a bialgebra section, implying that
$\Oc(\tilde{\Ga})\simeq \Oc(\Ga)\ot \CC[\partial N]$.

\smallbreak Let $\Lambda = \frac{1}{\vert \Upsilon \vert}\sum_{z \in
N} \delta(z)\partial^{-z}$ be the integral of $\CC[\Upsilon]$ and
denote by $L_{\Lambda}$ the endomorphism of $\Oc(\tilde{\Ga})$ given
by left multiplication of $\Lambda$. Since $\Oc(\tilde{\Ga})\simeq
\Oc(\Ga)\ot \CC[\partial N] \simeq \Oc(\Ga)\ot \CC[\Upsilon]$, it
follows that $\Ker L_{\Lambda} = \Oc(\Ga)(\CC[\Upsilon])^{+}$. But
since $A_{\lgot,\sigma} = \Oc(\tilde{\Ga}) \oplus M$ as
$\Oc(\tilde{\Ga})$-modules, we have that $J_{\delta} \cap
\Oc(\tilde{\Ga}) = A_{\lgot,\sigma}(\CC[\Upsilon])^{+} \cap
\Oc(\tilde{\Ga}) = \Oc(\tilde{\Ga})(\CC[\Upsilon])^{+} =
\Oc(\Ga)(\CC[\Upsilon])^{+} = \Ker L_{\Lambda}$. Hence $J_{\delta}
\cap \Oc(\Ga) = \Ker L_{\Lambda} \cap \Oc(\Ga) = 0$ for if $x \in
\Ker L_{\Lambda} \cap \Oc(\Ga)$, then $$ 0  = \Lambda x =
\frac{1}{\vert \Upsilon \vert}\sum_{z \in N} (\delta(z)\ot
\partial^{-z})(x\ot 1) = \frac{1}{\vert \Upsilon \vert}\sum_{z \in
N} \delta(z)x\ot
\partial^{-z},
$$

\noindent which implies that $\delta(z)x = 0$ for all $z\in N$,
because the elements $\partial^{z}$ are linearly independent. Thus
$x=0$ since $\delta(z)$ is invertible for all $z\in N$. \epf

\begin{obs} $(a)$
If $\Ga$ is finite-dimensional, then $\Oc(\Ga) = \CC^{\Ga}$ and by
Remark \ref{rmk:kt}, $\dim A_{\D}= \vert \Ga \vert \dim H$. In this
case, $\D$ is a finite subgroup datum and the last step of the proof
of the theorem above follows easily by dimension arguments. Indeed,
by \cite[Lemma 4.8]{Mu2}, we have that $\dim A_{\D} = \dim A_{\lgot,
\sigma} / \vert \Upsilon\vert$. Since $A_{\lgot, \sigma}$ and
$A_{\D}$ are extensions, it follows that
\begin{equation}\label{ec:dimensiones}\dim \CC^{\Ga}\frac{\dim
\ql}{\vert \Upsilon\vert} = \dim A_{\D} = \dim (\CC^{\Ga} / \J )\dim
H = \dim (\CC^{\Ga} / \J) \frac{\dim \ql}{\vert
N\vert}.\end{equation}

\noindent Since $\bar{\pi}(\Upsilon) = \{D^{z}\vert\ z\in N\}$ and
$\bar{\pi}(\partial^{z}\delta(z)^{-1}) = D^{z} = 1$ if and only if
$z = 0$, we have that $\vert \Upsilon \vert = \vert N \vert$. Thus,
from the equality \eqref{ec:dimensiones} it follows that $\CC^{\Ga}
= \CC^{\Ga}/\J$.

\smallbreak $(b)$ All exact sequences in the rows of diagram
\eqref{diag:const} are of the type $\mathcal{B}\hookrightarrow
\mathcal{A} \twoheadrightarrow \mathcal{H}$, where $\mathcal{B}$
is central in $\mathcal{A}$ and $\mathcal{H}$ is
finite-dimensional. Thus, by \cite[Thm. 1.7]{KT}, $\mathcal{B}
\subset \mathcal{A}$ is an $\mathcal{H}$-Galois extension and
$\mathcal{A}$ is a finitely-generated projective
$\mathcal{B}$-module. Moreover, using Lemma \ref{lema:quotient}
and Proposition \ref{prop:sucUl} $(b)$, one can see that the first
three exact sequences are cleft.
\end{obs}


\subsection{Relations between quantum
subgroups}\label{sec:equivalence}

Let $U$ be any Hopf algebra and consider the category $\quot (U)$,
whose objects are surjective Hopf algebra maps $q: U \to A$. If $q:
U \to A$ and $q': U \to A'$ are such maps, then an arrow
$\xymatrix{q\ar[0,1]^{\alpha}& q'}$ in $\quot (U)$ is a Hopf algebra
map $\alpha: A\to A'$ such that $\alpha q = q'$. In this language, a
\emph{quotient} of $U$ is just an isomorphism class of objects in
$\quot (U)$; let $[q]$ denote the class of the map $q$. There is a
partial order in the set of quotients of $U$, given by $[q]\leq
[q']$ iff there exists an arrow $\xymatrix{q\ar[0,1]^{\alpha}& q'}$
in $\quot (U)$. Notice that $[q]\leq [q']$ and $[q']\leq [q]$
implies $[q] = [q']$.

Our aim is to describe the partial order in the set $[q_{\D}]$,
$\D$ a subgroup datum, of quotients $q_{\D}: \Oe\twoheadrightarrow
A_{\D}$ given by Theorem \ref{teo:constrfinal}. Eventually, this
will be the partial order in the set of all quotients of $\Oe$. We
begin by the following definition. By an abuse of notation we
write $[A_{\D}] = [q_{\D}]$.

\begin{definition}\label{defi:order-datum} Let
$\D= (I_{+}, I_{-}, N, \Ga, \sigma, \delta)$ and $\D'= (I'_{+},
I'_{-}, N', \Ga', \sigma', \delta')$ be subgroup data. We say that
$\D \leq \D'$ iff

\begin{itemize}
\item[$\bullet$] $I'_{+} \subseteq I_{+}$ and $I'_{-} \subseteq
I_{-}$.

\noindent In particular, this condition implies that $I'\subseteq
I$, $\mathbb{T}_{I'} \subseteq \mathbb{T}_{I}$ and
$\mathbb{T}_{I^{c}} \subseteq \mathbb{T}_{I'^{c}}$. Since $\Sigma
= \mathbb{T}_{I}\times \Omega$ and $\Sigma' =
\mathbb{T}_{I'}\times \Omega'$, we have that $\Omega'\subseteq
\Omega \subseteq \mathbb{T}_{I^{c}} \subseteq
\mathbb{T}_{I'^{c}}$. As $\mathbb{T}_{I'^{c}} = \mathbb{T}_{I^{c}}
\times \mathbb{T}_{I'^{c}-I^{c}}$, the restriction map
$\widehat{\mathbb{T}_{I'^{c}}} \twoheadrightarrow
\widehat{\mathbb{T}_{I^{c}}}$ admits a canonical section $\eta$
and $\eta(N) \subseteq N'$.

\smallbreak\item[$\bullet$] There exists a morphism of algebraic
groups $\tau:\Ga' \to \Ga$ such that $\sigma \tau = \sigma'$.

\smallbreak\item[$\bullet$] $\delta'\eta = \ ^{t}\tau \delta$.

\end{itemize}
Furthermore, we say that $\D \simeq \D'$ iff $\D \leq \D'$ and $\D'
\leq \D$. This means that

\begin{itemize}
\item[$\bullet$] $I_{+} = I'_{+}$ and $I_{-} = I'_{-}$.

\smallbreak \item[$\bullet$] There exists an isomorphism of
algebraic groups $\tau:\Ga' \to \Ga$ such that $\sigma \tau =
\sigma'$.

\smallbreak\item[$\bullet$] $N = N'$ and $\delta' = \ ^{t}\tau
\delta$.
\end{itemize}
\end{definition}

\begin{theorem}\label{teo:order-datum} Let
$\D$ and $\D'$ be subgroup data. Then

(a) $[A_{\D}] \leq [A_{\D'}]$ iff $\D\leq\D'$.

(b) $[A_{\D}] = [A_{\D'}]$ iff $\D\simeq\D'$.
\end{theorem}
\pf Let $q = q_{\D}$ and $q' = q_{\D'}$. Suppose that $[A_{\D}]
\leq [A_{\D'}]$, that is, there exists a surjective Hopf algebra
map $\alpha : A_{\D} \to A_{\D'}$ such that $\alpha q = q'$. Since
by Theorem \ref{teo:constrfinal}, $\hat{\iota}\ ^{t}\sigma =
q\iota$ and $\hat{\iota}'\ ^{t}\sigma' = q'\iota$, we have that
$\alpha \hat{\iota}\ ^{t}\sigma $ $= \alpha q\iota = q' \iota =
\hat{\iota}'\ ^{t}\sigma'$. Thus, the Hopf algebra map $\beta:=
\alpha \hat{\iota} : \Oc(\Ga) \to \Oc(\Ga')$ is surjective with
$\Img \beta \subseteq \Img \ ^{t}\sigma$ and its transpose defines
an injective map of algebraic groups $\tau : \Ga' \to \Ga$ such
that $\sigma\tau=\sigma'$.

\smallbreak Again by Theorem \ref{teo:constrfinal}, we know that
both $A_{\D}$ and $A_{\D'}$ are central extensions by  $H\simeq
A_{\D} /A_{\D}\Oc(\Ga)^{+}$ and $H'\simeq A_{\D'}
/A_{\D'}\Oc(\Ga')^{+}$, respectively. Since $\hat{\pi}' \alpha
(A_{\D}\Oc(\Ga)^{+}) = \hat{\pi}' (A_{\D'}\Oc(\Ga')^{+}) = 0$, there
exists a surjective Hopf algebra map $\gamma: H \to H'$ such that
the following diagram commutes
$$\xymatrix{1 \ar[r]^{} & \Oc(G)
\ar[r]^{\iota}\ar@/_2pc/@{.>}[dd]_(.7){^{t}\sigma'}
\ar[d]_{^{t}\sigma} & \Oe \ar[r]^{\pi}
\ar[d]^{q}\ar@/_2pc/@{.>}[dd]_(.7){q'} & \qe^{*}
\ar[r]^{}\ar[d]_{r}
\ar@/^2pc/@{.>}[dd]^(.7){r'} & 1\\
1 \ar[r]^{} &   \Oc(\Ga)\ar[d]^{\beta}  \ar[r]^{\hat{\iota}} &
A_{\D} \ar[r]^{\hat{\pi}}\ar[d]^{\alpha} &
H \ar[r]^{}\ar[d]_{\gamma} & 1\\
1 \ar[r]^{} &   \Oc(\Ga')  \ar[r]^{\hat{\iota}'} &
A_{\D'}\ar[r]^{\hat{\pi}'} & H' \ar[r]^{} & 1.}$$

\noindent Since $^{t}r: H^{*} \hookrightarrow \qe$ and
$^{t}r':(H')^{*} \hookrightarrow \qe$ are just the inclusions, it
follows that $^{t}\gamma : (H')^{*} \hookrightarrow H^{*}$ is the
same inclusion. If $H^{*}$ and $(H')^{*}$ are determined by the
triples $(\Sigma, I_{+}, I_{-})$ and $(\Sigma', I_{+}', I_{-}')$,
it follows that $\Sigma'\subseteq \Sigma$, $I_{+}' \subseteq
I_{+}$, $I_{-}' \subseteq I_{-}$, whence $\eta(N) \subseteq N'$.
Thus, $\qlp \subseteq \ql$ by Lemma \ref{lema:uelsualg}.

\smallbreak Now by Theorem \ref{teo:constrfinal}, $\delta(z) =
t(\partial^{z})$ in $A_{\D}$ and $\delta'(z') = t'(\partial^{z'})$
in $A_{\D'}$, for all $z \in N$ and $z'\in N'$. Thus, for all $z \in
N$ we have
\begin{align*}
^{t}\tau \delta(z) = \alpha \delta (z) = \alpha t(\partial^{z}) =
\alpha t\nu(\psi^{*}(D^{z})) = t'\nu'((\psi')^{*}\eta(D^{z})) =
\delta'(\eta(z)),
\end{align*}

\noindent where the fourth equality follows from the construction of
the quotients $A_{\D}$, $A_{\D'}$ and $\alpha q = q'$. All this
implies that $\D\leq \D'$.

\smallbreak Suppose now that $\D\leq \D'$. This implies that $\qlp
\subseteq \ql$ and by construction, there exists a Hopf algebra
map $\kappa: \Ol \to \Olp$ such that
$$\xymatrix{\Oe \ar@{->>}[r]^{\Res}\ar@{->>}[dr]_{\Res'} &
\Ol\ar@{->>}[d]^{\kappa}\\
& \Olp}$$

\noindent commutes. Since  $ ^{t}\tau\ ^{t}\sigma=\ ^{t}\sigma'$,
there exists a commutative diagram
$$\xymatrix{\Oc(L)
\ar[rr]^{\iota_{L}} \ar[d]_{^{t}\sigma} & &
\Ol \ar[d]_{\nu} \ar@/^/[ddr]^{t'\nu'\kappa}\\
 \Oc(\Ga)\ar[rr]_{\bar{\iota}} \ar@/_/[dr]_{^{t}\tau}& & A_{\lgot,\sigma}
\\
& \Oc(\Ga')\ar[rr]_{\hat{\iota}'} & & A_{\D'}.}$$

\noindent As $A_{\lgot,\sigma}$ is a pushout, there exists a
surjective Hopf algebra map $\tilde{\alpha} :  A_{\lgot,\sigma} \to
A_{\D'} $ such that $\tilde{\alpha} \nu = t'\nu'\kappa$. Since
$A_{\D} = A_{\lgot,\sigma} / J_{\delta}$, to show the existence of a
surjective map $\alpha: A_{\D} \to A_{\D'}$ such that $\alpha q =
q'$, it is enough to prove that $\tilde{\alpha} (J_{\delta}) = 0$.
But $J_{\delta}$ is the two-sided ideal of $A_{\lgot,\sigma}$
generated by $\delta(z) - \partial^{z}$ with $z\in N$; now
\begin{align*}
\tilde{\alpha}(\delta(z) -
\partial^{z})  & =\ ^{t}\tau\delta(z) -
\tilde{\alpha}(\nu\psi^{*}(D^{z})) = \ ^{t}\tau\delta(z) -
t'\nu'\eta(z) \\
&=  \ ^{t}\tau\delta(z) - \delta'\eta(z) = 0,
\end{align*}
\noindent  by assumption. Hence, $\tilde{\alpha} (J_{\delta}) =
0$. This finishes the proof of $(a)$. Now $(b)$ follows
immediately. \epf



\section{Determining quantum subgroups}\label{sec:demmainthm}
Let $q: \Oe \to A$ be a surjective Hopf algebra map. We prove now
that it is isomorphic to $q_{\D}: \Oe \to A_{\D}$ for some subgroup
datum $\D$. This concludes the proof of Theorem \ref{teo:biyeccion}.

\smallbreak The Hopf subalgebra $K= q(\Oc(G))$ is central in $A$ and
whence $A$ is an $H$-extension of $K$, where $H$ is the Hopf algebra
$H = A/ AK^{+}$. Indeed, it follows directly from \cite[Prop.
3.4.3]{Mo}, because $A$ is faithfully flat over $K$ by \cite[Thm.
3.3]{Sch2}. Since $K$ is a quotient of $\Oc(G)$, there exists an
algebraic group $\Ga$ and an injective map of algebraic groups
$\sigma: \Ga \to G$ such that $K \simeq \Oc(\Ga)$. Moreover, since
$q(\Oe\Oc(G)^{+}) = AK^{+}$, we have that $\Oe\Oc(G)^{+} \subseteq
\Ker \hat{\pi}q$, where $\hat{\pi}: A \to H$ is the canonical
projection. Since $\qe^{*} \simeq \Oe /[\Oe\Oc(G)^{+}]$, there
exists a surjective map $r:\qe^{*} \to H$ and by Proposition
\ref{prop:hopf-subalg}, $H^{*}$ is determined by a triple $(\Sigma,
I_+, I_-)$. In particular, we have the following commutative diagram
\begin{equation}\label{diag:propqsubgr}
 \xymatrix{1
\ar[r]^{} & \Oc(G) \ar[r]^{\iota} \ar[d]_{^{t}\sigma}  & \Oe
\ar[r]^{\pi}
\ar[d]^{q} & \qe^{*} \ar[r]^{}\ar[d]^{r} & 1\\
1 \ar[r]^{} &   \Oc(\Ga) \ar[r]^{\hat{\iota}} & A
\ar[r]^{\hat{\pi}} & H \ar[r]^{} & 1.}
\end{equation}

\smallbreak Let $N$ correspond to $\Sigma$ as in Remark
\ref{rmk:subgruposigma}. Our aim is to show that there exists
$\delta$ such that $A \simeq A_{\D}$ for the subgroup datum $\D =
(I_{+}, I_{-}, N, \Ga, \sigma, \delta)$. Recall the Lie algebra
$\liel$ from Definition \ref{def:subgroupdatum} and the Hopf
algebra $\ql\supseteq H^*$ from  \ref{subsub:ul}. Denote by $v :
\ql^{*} \to H$ the surjective Hopf algebra map induced by this
inclusion.

\begin{lema}\label{prop:fact}
The diagram \eqref{diag:propqsubgr} factorizes through the exact
sequence
$$\xymatrix{1
\ar[r]^{} & \Oc(L) \ar[r]^{\iota_{L}} & \Oel \ar[r]^{\pi_{L}} &
\ql^{*} \ar[r]^{}& 1,}$$

\noindent that is, there exist Hopf algebra maps $u,\ w$ such that
the following diagram with exact rows commutes:

$$\xymatrix{1
\ar[r]^{} & \Oc(G)
\ar[r]^{\iota}\ar@/_2pc/@{.>}[dd]_(.7){^{t}\sigma} \ar[d]_{\res} &
\Oe \ar[r]^{\pi} \ar[d]^{\Res}\ar@/_2pc/@{.>}[dd]_(.7){q} & \qe^{*}
\ar[r]^{}\ar[d]_{p}
\ar@/^2pc/@{.>}[dd]^(.7){r} & 1\\
1 \ar[r]^{} &   \Oc(L)\ar[d]^{u}  \ar[r]^{\iota_L} & \Ol
\ar[r]^{\pi_L}\ar[d]^{w} &
\ql^* \ar[r]^{}\ar[d]_{v} & 1\\
1 \ar[r]^{} &   \Oc(\Ga)  \ar[r]^{\hat{\iota}} &
A\ar[r]^{\hat{\pi}} & H \ar[r]^{} & 1.}$$
\end{lema}

\pf To show the existence of the maps $u$ and $w$ it is enough to
show that $\Ker \Res \subseteq \Ker q$, since $u$ is simply $w
\iota_{L}$. This clearly implies that $v\pi_L = \hat\pi w$.

\smallbreak Let $\QEbmas$ and $\QEbmen$ be the Borel subalgebras of
$\QEpe$ (see \cite{DL} and \cite[Cap. 4]{J}), and let
$\AA_{\epsilon}$ be the subalgebra of $\QEbmas \ot \QEbmen$
generated by the elements
$$\{1\ot e_{j}, f_{j}\ot 1, K_{-\lambda}\ot K_{\lambda}:\ 1\leq
j\leq n, \lambda \in P \},$$

\noindent where $P$ is the weight lattice. By \cite[Sec. 4.3]{DL},
this algebra has a basis given by the set $\{fK_{-\lambda}\ot
K_{\lambda}e \}$, where $\lambda \in P$ and $e$, $f$ are monomials
in $e_{\alpha}$ and $f_{\beta}$ respectively, $\alpha,\ \beta \in
Q_{+}$. Moreover, $\AA_{\epsilon}$ is a $(Q_{-},P,Q_{+})$-graded
algebra whose gradation is given by
\begin{align*}
\deg (f_{j} \ot 1 )  = (-\alpha_{j},& 0,0), \qquad \qquad \deg (1
\ot
e_{j} ) = (0,0,\alpha_{j}),\\
&\deg (K_{-\lambda}\ot K_{\lambda}) =
(0,\lambda,0),
\end{align*}

\noindent for all $1\leq j\leq n$, $\lambda \in P$. By \cite[4.3 and
6.5]{DL}, there exists an injective algebra map $\mu_{\epsilon} :
\Oe \to \AA_{\epsilon}$ such that $\mu_{\epsilon}(\Oc(G)) \subseteq
\AA_{0}$, where $\AA_{0}$ is the subalgebra of $\AA_{\epsilon}$
generated by the elements
$$\{1\ot e_{j}^{\ell}, f_{j}^{\ell}\ot 1, K_{-\ell\lambda}\ot
K_{\ell\lambda}:\ 1\leq j\leq n, \lambda \in P \}.$$

\noindent Hence, it is enough to show that $\mu_{\epsilon}(\Ker
\Res) \subseteq \mu_{\epsilon}(\Ker q)$.

\smallbreak \noindent {\bf Claim:} $\mu_{\epsilon}(\Ker \Res)$ is
the two-sided ideal $\II$ generated by the elements $$\{1\ot e_{k},
f_{j}\ot 1:\ \alpha_{k}\notin I_{-}, \alpha_{j}\notin I_{+}\}.$$

\noindent Indeed, let $\lambda \in P_{+}$ and let $\psi_{\lambda}
\in \Glie^{\circ}$ such that
\begin{align*}
\psi_{\lambda}^{}(FME) & =
\delta_{1,E}\delta_{1,F}M(\lambda),\qquad &
\qquad\psi_{-\lambda}^{}(EMF) & =
\delta_{1,E}\delta_{1,F}M(-\lambda),
\end{align*}
\noindent for all elements $FME$ of the PBW basis of $\Glie$, where
$M \in Q$ and the form $M(\lambda)$ is simply the linear extension
of the bilinear form  $<\alpha_{j}, \lambda> =
\epsilon^{d_{i}(\alpha_{i}, \lambda)}$ for all $\lambda \in P$,
$1\leq i\leq n$. By \cite[Sec. 4.4]{DL}, there exist matrix
coefficients $\psi_{\pm\lambda}^{\pm\alpha}$, and $\alpha \in Q_{+}$
such that
\begin{align*}
\psi_{-\lambda}^{\alpha}(EMF) & =
\psi_{-\lambda}^{}(EMFE_{\alpha}),\qquad &
\qquad\psi_{-\lambda}^{-\alpha}(EMF) & =
\psi_{-\lambda}^{}(F_{\alpha}EMF),
\end{align*}

\noindent for all elements $EMF$ of the PBW basis of $\Glie$.
Moreover, one has that
\begin{align*}
\mu_{\epsilon}(\psi_{-\varpi_{i}})  = K_{-\varpi_{i}}\ot
&K_{\varpi_{i}}, \qquad\qquad
\mu_{\epsilon}(\psi_{-\varpi_{i}}^{\alpha_{k}}) =
K_{-\varpi_{i}}\ot K_{\varpi_{i}}e_{k},\\
& \mu_{\epsilon}(\psi_{-\varpi_{i}}^{-\alpha_{j}}) =
f_{j}K_{-\varpi_{i}}\ot K_{\varpi_{i}},
\end{align*}

\noindent for all $1\leq i,j\leq n$. Through a direct computation
one can see that $\psi_{-\varpi_{i}}^{\alpha_{k}},\
\psi_{-\varpi_{i}}^{-\alpha_{j}} \in \Ker \Res$ and

$$\mu_{\epsilon}(\psi_{\varpi_{i}}\psi_{-\varpi_{i}}^{\alpha_{k}}) =
1\ot e_{k}\qquad
\mu_{\epsilon}(\psi_{-\varpi_{i}}^{-\alpha_{j}}\psi_{\varpi_{i}}) =
f_{j}\ot 1.$$

\noindent for all $\alpha_{k} \notin I_{-}$, $\alpha_{j}\notin
I_{+}$. Hence, the generators of $\II$ are in $\mu_{\epsilon}(\Ker
\Res)$.

\smallbreak Conversely, if $h \in \Ker \Res$, then $h|_{\Glil} =0$
and by definition we have that
$$<\mu_{\epsilon}(h), EM\ot NF> = <h, EMNF> = 0,$$

\noindent for all elements $EMNF$ of the PBW basis of $\Glil$. Thus,
using the existence of perfect pairings (see \cite[Sec. 3.2]{DL})
and evaluating in adequate elements, it follows that each term of
the basis $\{fK_{-\lambda}\ot K_{\lambda}e \}$ that appears in
$\mu_{\epsilon}(h)$ must lie in $\II$.

\smallbreak Since $0 = \pi_{L}\Res (h) = r\pi (h) = \hat{\pi}q(h)$,
we have that $q(h) \in \Ker \hat{\pi} = \Oc(\Ga)^{+}A =
q(\Oc(G)^{+}\Oe)$. Then there exist $a\in \Oc(G)^{+}\Oe$ and $c \in
\Ker q$ such that $h = a + c$; in particular, for all generators $t$
of $\II$ we have that $ t = \mu_{\epsilon}(a) + \mu_{\epsilon}(c)$,
where $\mu_{\epsilon}(a)$ is contained in $\AA_{0}$. Comparing
degrees in both sides of the equality we have that
$\mu_{\epsilon}(a) = 0$, which implies that each generator of $\II$
must lie in $\mu_{\epsilon}(\Ker q)$. \epf

The following lemma shows the convenience of characterizing the
quotients $A_{\lgot, \sigma}$ of $\Oe$ as pushouts.

\begin{lema}\label{lema:Acociente}
$\sigma(\Ga)\subseteq L$ and therefore $A$ is a quotient of
$A_{\lgot, \sigma}$ given by the pushout. Moreover, the following
diagram commutes
\begin{equation}\label{diag:gordo}
\xymatrix{1 \ar[r]^{} & \Oc(G) \ar[r]^{\iota} \ar[d]_{\res} & \Oe
\ar[r]^{\pi}
\ar[d]^{\Res} & \qe^{*} \ar[r]^{}\ar[d]^{p} & 1\\
1 \ar[r]^{} &   \Oc(L)\ar[d]_{u}  \ar[r]^{\iota_L} & \Ol
\ar[r]^{\pi_L}\ar[d]^{\nu} &
\ql^* \ar[r]^{}\ar@{=}[d] & 1\\
1 \ar[r]^{} &   \Oc(\Ga)  \ar[r]^{j}\ar@{=}[d] & A_{\lgot,
\sigma}\ar[r]^{\bar{\pi}}\ar[d]^{t} & \ql^* \ar[r]^{}\ar[d]^{v} & 1\\
1 \ar[r]^{} &   \Oc(\Ga)  \ar[r]^{\hat{\iota}} &
A\ar[r]^{\hat{\pi}}& H\ar[r]^{}& 1.}
\end{equation}
\end{lema}

\pf Recall the maps $u,\ w$ defined in the lemma above; we have
that $w\iota_{L} = \hat{\iota}u$, that is, the following diagram
commutes

$$\xymatrix{\Oc(L)
\ar[r]^{\iota_{L}} \ar[d]_{u} &
\Ol \ar[d]^{\nu} \ar@/^/[ddr]^{w}\\
\Oc(\Ga) \ar[r]_{j} \ar@/_/[drr]_{\hat{\iota}}& A_{\lgot, \sigma}
\\
& & A.}$$

\noindent Since $A_{\lgot, \sigma}$ is a pushout, there exists a
unique Hopf algebra map $t: A_{\lgot, \sigma} \to A$ such that $ts
= w$ and $tj = \hat{\iota}$. This implies that $\Ker \bar{\pi} =
j(\Oc(\Ga))^{+}A_{\lgot,\sigma} \subseteq \Ker \hat{\pi}t$ and
therefore the diagram \eqref{diag:gordo}  is commutative. \epf

Let $(\Sigma, I_{+}, I_{-})$ be the triple that determines $H$.
Recall that by Remark \ref{rmk:subgruposigma}, giving a group
$\Sigma$ such that $\mathbb{T}_{I} \subseteq \Sigma \subseteq
\mathbb{T}$ is the same as giving a subgroup $N\subseteq
\widehat{\mathbb{T}_{I^{c}}}$. In fact, by Lemma \ref{lema:grupoz},
we know that the Hopf algebra $A_{\lgot, \sigma}$ contains a set of
central group-like elements $\z= \{\partial^{z}\vert\ z\in
\widehat{\mathbb{T}_{I^{c}}}\}$ such that $\bar{\pi}(\partial^{z}) =
D^{z}$ for all $z\in \widehat{\mathbb{T}_{I^{c}}}$ and $H =
\ql^{*}/(D^{z} -1\vert\ z\in N)$. To see that $A = A_{\D}$ for a
subgroup datum $\D = (I_{+}, I_{-}, N, \Ga, \sigma, \delta)$ it
remains to find a group map $\delta: N \to \widehat{\Ga}$ such that
$A \simeq A_{\lgot, \sigma}/ J_{\delta}$. This is given by the last
lemma of the paper.

\begin{lema}\label{lema:thetaiso}
There exists a group homomorphism $\delta: N \to \widehat{\Ga}$ such
that $J_{\delta} = (\partial^{z} - \delta(z)\vert\ z\in N)$ is a
Hopf ideal of $A_{\lgot, \sigma}$ and $A\simeq A_{\D}= A_{\lgot,
\sigma}/ J_{\delta}$.
\end{lema}

\pf Let $\partial^{z} \in \z$. Then $\hat{\pi}t(\partial^{z}) =
v\bar{\pi}(\partial^{z}) = 1$ for all $z \in N$, by Lemma
\ref{lema:elemcentrales} $(b)$. Since $t(\partial^{z})$ is a
group-like element, this implies that $t(\partial^{z}) \in A^{\co
\hat{\pi}} = \Oc(\Ga)$. As $G(\Oc(\Ga)) = \widehat{\Ga}$, we have a
group homomorphism $\delta$ given by
$$\delta: N \to \widehat{\Ga},\qquad \delta(z) =
t(\partial^{z})\quad\forall\ z\in N.$$

The two-sided ideal of $A_{\lgot,\sigma}$ given by $J_{\delta} =
(\partial^{z}-\delta(z)\vert\ z\in N)$ is clearly a Hopf ideal and
$t(J_{\delta}) = 0$. Consequently we have a surjective Hopf
algebra map $ \theta: A_{\D}  \twoheadrightarrow A$, which makes
the following diagram commutative
\begin{equation}\label{diag:episeq}
\xymatrix{1 \ar[r]^{} &   \Oc(\Ga)
\ar[r]^{\tilde{\iota}}\ar@{=}[d] & A_{\D}
\ar[r]^{\tilde{\pi}}\ar@{->>}[d]^{\theta} & H\ar[r]^{}\ar@{=}[d] & 1\\
1 \ar[r]^{} &  \Oc(\Ga)  \ar[r]^{\hat{\iota}} & A
\ar[r]^{\hat{\pi}}& H\ar[r]^{}& 1.}
\end{equation}
Then $\theta$ is an isomorphism by Corollary \ref{cor:five}. \epf

\subsection*{Acknowledgments} We thank Akira Masouka
for kindly communicating us Lemma
\ref{lema:iso-galois-extensions}.



\bibliographystyle{amsbeta}

\end{document}